\documentclass[11pt]{article}
\usepackage{graphicx,wrapfig}
\usepackage{flafter}
\usepackage{leftidx}
\usepackage{mathrsfs}
\usepackage{amssymb,amsmath}
\usepackage{bbding}
\usepackage{fancyhdr}
\usepackage{mathrsfs}
\usepackage{color}
\usepackage{stmaryrd}
\usepackage{amsfonts}
\usepackage{latexsym}
\usepackage{psfrag}
\usepackage{graphicx,subfigure}
\usepackage{fancyhdr,graphicx}
\usepackage{multicol}
\usepackage{dsfont}
\usepackage{bbm}
\usepackage{booktabs}
\usepackage[center]{caption2}
\usepackage{cite}
\usepackage{epstopdf}
\usepackage{booktabs}
\usepackage{multirow}
\usepackage{enumerate}
\usepackage{enumitem}
\usepackage{algpseudocode,algorithm,algorithmicx}

\allowdisplaybreaks

\date{}
\newcommand{\zd}{\,\mathrm{d}}

\newtheorem{theorem}{Theorem}[section]
\newtheorem{lemma}{Lemma}[section]

\newtheorem{example}{Example}[section]
\newtheorem{corollary}{Corollary}[section]
\numberwithin{equation}{section}

\newcommand{\diff}{\triangledown_\tau}
\newcommand{\abs}[1]{\left|#1\right|}
\newcommand{\absb}[1]{\big|#1\big|}

\newcommand{\bra}[1]{\left(#1\right)}
\newcommand{\brab}[1]{\big(#1\big)}
\newcommand{\braB}[1]{\Big(#1\Big)}

\newcommand{\kbra}[1]{\left[#1\right]}
\newcommand{\kbrab}[1]{\big[#1\big]}

\newcommand{\mynorm}[1]{\left\|#1\right\|}
\newcommand{\mynormb}[1]{\big\|#1\big\|}
\newcommand{\mynormB}[1]{\Big\|#1\Big\|}

\begin{document}
\title{A second-order and nonuniform time-stepping maximum-principle preserving scheme for time-fractional Allen-Cahn equations}
\author{Hong-lin Liao\thanks{
Department of Mathematics,
Nanjing University of Aeronautics and Astronautics,
Nanjing 211106, P. R. China. E-mails: liaohl@nuaa.edu.cn and liaohl@csrc.ac.cn.
Research supported by NUAA Scientific Research Starting Fund of Introduced Talent (1008-56SYAH18037).}
\quad Tao Tang\thanks{Department of Mathematics and International Center for Mathematics, Southern
    University of Science and Technology, Shenzhen, Guangdong Province; and
    Division of Science and Technology, BNU-HKBU United International College,
    Zhuhai, Guangdong Province, China.
    Email: tangt@sustech.edu.cn. This author's work is partially supported by the NSF of China under grant number 11731006.}
\quad Tao Zhou\thanks{NCMIS \& LSEC, Institute of Computational Mathematics and Scientific/Engineering Computing,
Academy of Mathematics and Systems Science, Chinese Academy of Sciences, Beijing, 100190,
P. R. China. Email: tzhou@lsec.cc.ac.cn. This author's work is partially supported by the NSF of China (under grant numbers 11822111, 11688101, 91630203, 11571351, and 11731006), the science challenge project (No. TZ2018001), NCMIS, and the youth innovation promotion association (CAS).}}
\date{}
\maketitle
\normalsize

\begin{abstract}
In this work, we present a second-order nonuniform time-stepping scheme for the time-fractional Allen-Cahn equation. We show that the proposed scheme preserves the discrete maximum principle, and by using the convolution structure of consistency error, we present sharp maximum-norm error estimates which reflect the temporal regularity. As our analysis is built on nonuniform time steps, we may resolve the intrinsic initial singularity by using the graded meshes. Moreover, we propose an adaptive time-stepping strategy for large time simulations. Numerical experiments are presented to show the effectiveness of the proposed scheme. This seems to be the first second-order maximum principle preserving scheme for the time-fractional Allen-Cahn equation.
\\

{\textsc{Keywords:}}\;\; Time-fractional Allen-Cahn equation;
Alikhanov formula; adaptive time-stepping strategy;
discrete maximum principle; sharp error estimate \\

\noindent{\bf AMS subject classiffications.}\;\; 35Q99, 65M06, 65M12, 74A50
\end{abstract}
\section{Introduction}
We consider the following two dimensional time-fractional Allen-Cahn equation
\begin{align}
\partial_{t}^{\alpha}u(\mathbf{x},0)
=&\,\varepsilon^{2}\Delta{u}-f(u),\quad
\mathbf{x}\in\Omega,\quad{0}<{t}\leq{T},\label{Problem-1}\\
u(\mathbf{x},0)=&\,u_{0}(\mathbf{x}),\quad \mathbf{x}\in\bar{\Omega},\label{Problem-2}
\end{align}
where $\Omega=(0,L)^2$ with closure $\bar{\Omega}$.  The nonlinear bulk force $f(u)$ is given by $f(u)=u^{3}-u.$ The constant $0<  \varepsilon  \ll 1$ is the interaction length that describes the thickness of the transition boundary between materials. For simplicity, we consider the periodic boundary conditions. In equations (1.1), $\partial_{t}^{\alpha}:={}_{0}^{C}\!D_{t}^{\alpha}$ denotes
the Caputo derivative of order $\alpha:$
\begin{align}\label{CaputoDef}
(\partial_{t}^{\alpha}v)(t)
:=(\mathcal{I}_{t}^{1-\alpha}v')(t)=\int_{0}^{t}\omega_{1-\alpha}(t-s)v'(s)\zd{s},\quad 0<\alpha<1,
\end{align}
where $\mathcal{I}_{t}^{\mu}$ is the fractional Riemann-Liouville integral of order $\mu>0$, that is,
\begin{align}
(\mathcal{I}_{t}^{\mu}v)(t)
:=\int_{0}^{t}\omega_{\mu}(t-s)v(s)\zd{s},\quad\text{where}\quad  \omega_{\mu}(t):=t^{\mu-1}/\Gamma(\mu).
\end{align}
As a generalization of the classical Allen-Cahn equation \cite{AllenCahn_1979,Guan.etal2014,ShenY_DCDS_2010,Feng_NumMath_2003}, the above time-fractional Allen-Cahn equation (1.1) has been widely investigated in recent years \cite{Inc2018Time,fractional_Wang_I,Liu2018Time,Tang2018On}, In particular, it was first shown in \cite{Tang2018On} that the time-fractional Allen-Cahn equation admits the following energy law
\begin{align}\label{Frac-Energy-DecayLaw}
E(t)\leq{E}(0),
\end{align}
where $E(t)$ is the total energy defined by
\begin{align}\label{Frac-Energy-Expression}
E(t):=\int_{\Omega}\left[\frac{\varepsilon^{2}}{2}|\nabla{u}|^{2}
+F(u)\right]\zd{\mathbf{x}},\quad
F(u)=\frac{1}{4}(1-u^{2})^{2}.
\end{align}
Moreover, the following maximum principle holds \cite{Tang2018On}
\begin{align}\label{Frac-Maximum-Principle}
|u(\mathbf{x},t)|\leq{1}\quad
\text{if}\quad|u(\mathbf{x},0)|\leq{1}.
\end{align}
From the numerical scheme point of view, first order schemes that combine the $L1$ formula \cite{lin_finite_2007,sun_fully_2006} and the stabilization technique \cite{XuT_SINUM_2006} were proposed in \cite{Tang2018On} for the time-fractional Allen-Cahn equation. Furthermore, it is shown in \cite{Tang2018On} that the stabilization $L1$ scheme preserves the energy law (\ref{Frac-Energy-DecayLaw}) and the maximum principle (\ref{Frac-Maximum-Principle}) in the discrete level. More recently, sharp regularity analysis of the time-fractional Allen-Cahn equation and numerical analysis for a class of numerical schemes under limited regularity were presented in \cite{Du}. Notice that all the analysis in the above mentioned works is based on uniform time grids.

In this work, we aim at designing a second order scheme using \textit{nonuniform} time grids. There are two main motivations to investigate nonuniform time grids:
\begin{itemize}
\item Similar as in other nonlinear subdiffusion problems, the solution of the time-fractional Allen-Cahn equation admits an intrinsic initial singularity \cite{Jin2016An}. Consequently, graded meshes are recommended for such problems \cite{stynes_error_2017,Martin}.
\item The solution of the time-fractional Allen-Cahn equation may admit multiple time scales \cite{Gomez2011Provably,QiaoAn2011,Li2017Computationally,Jimack,Xucicp}, i.e., the initial dynamics evolve on a fast time scale and later coarsening stage evolves on a very slow time scale. Therefore, one may need to use adaptive time grids to catch different time scales \cite{QiaoAn2011,Zhang2012An}.
\end{itemize}

To this end, we present in this work a second order Alikhanov-type scheme under nonuniform time grids. We shown that the proposed scheme preserves the discrete maximum principle, and this seems to be the first work on  second order maximum principle preserving schemes for the time-fractional Allen-Cahn equation. We also present a sharp maximum-norm error estimate the can reflect the temporal regularity. Finally, we propose an adaptive time-stepping strategy for long-time simulations.

The rest of this work is organized as follows. We present some preliminaries in Section 2.  A second order nonuniform Alikhanov scheme is proposed in Section 3, where the discrete maximum principle is also established.  The convergence analysis of the proposed scheme is given in section 4, and this is followed by
extensive experiments in Section 5. We finally give some concluding remarks in Section 6.

\section{Preliminaries}
In this section, we shall present some preliminaries.
\subsection{Nonuniform time grids}
Throughout the paper, we shall consider nonuniform time grids. To this end, we introduce the following time mesh:
\begin{equation}{\label{nonuniformgrid}}
0=t_{0}<t_{1}<\cdots<t_{k-1}<t_{k}<\cdots<t_{N}=T
\end{equation}
with time-step sizes $\tau_{k}:=t_{k}-t_{k-1}$ for $1\leq{k}\leq{N}.$ We define the maximum time-step size as $\tau:=\max_{1\leq{k}\leq{N}}\tau_{k}$.
Also, for $k\geq{1}$ and $0<\theta<1$ we define the off-set time level as $t_{k-\theta}:=(1-\theta)t_{k}+\theta t_{k-1}.$
We set the adjacent step ratio as $\rho_k:=\tau_k/\tau_{k+1}$ and define the maximum step ratio as $\rho:=\max_{k\geq 1}\rho_k$. We now introduce the following assumptions.
\begin{enumerate}[itemindent=1em]
\item[\textbf{M1}.] The maximum time-step ratio $\rho=7/4$.
\end{enumerate}
The condition \textbf{M1} says that one can use a series of decreasing time-steps with the reduction factor down to $4/7$.
Always, we do not impose any restrictions to the amplification factor for increasing time-steps,
although a maximum time-step size may be necessary for theocratical analysis. The use of nonuniform meshes are motivated by the following two reasons:
Firstly, to resolve the initial solution singularity $u_{t} \sim \mathcal{O}(t^{\alpha-1})$ as $t\rightarrow0$,
a graded mesh $t_k=T(k/N)^{\gamma}$ with the step ratios $\rho_k\le1$
has been a  popular approach in recent years \cite{stynes_error_2017}. Secondly, to capture the fast dynamics away from $t=0$
and the slowly coarsening stage near the steady state, one may use an adaptive time-stepping strategy \cite{QiaoAn2011,Li2017Computationally}.
We also need the following assumption for the sake of convergence analysis \cite{Liao2018Sharp,Liao2018Unconditional,William2007A}:
\begin{enumerate}[itemindent=1em]
\item[\textbf{M2}.] For a parameter $\gamma\geq{1}$, there exists mesh-independent constants $C_{1\gamma},C_{2\gamma}>0$ such that
$\tau_k\le \tau\min\{1,C_{1\gamma}t_k^{1-1/\gamma}\}$
for~$1\le k\le N$~and
$t_{k}\leq C_{2\gamma}t_{k-1}$ for $2\leq{k}\leq{N}$.
\end{enumerate}
Here, the parameter $\gamma\ge1$ controls the extent to which the time levels
are concentrated near $t=0$. If the mesh is quasi-uniform, then
\textbf{M2} holds with~$\gamma=1$.  As~$\gamma$ increases, the initial step sizes
become smaller compared to the later ones.

To facilitate the error analysis of difference approximations in space, we assume that the continuous solution $u$
is sufficiently smooth in physical domain and satisfies
\begin{align}\label{Regularity-Sigma}
\mynormb{u(t)}_{W^{4,\infty}(\Omega)}\leq{C}_u,\quad
\mynormb{u^{(\ell)}(t)}_{W^{2,\infty}(\Omega)}\leq{C}_u\brab{1+t^{\sigma-\ell}}\;\;
\text{for $\ell=1,2,3$,}
\end{align}
where a regularity parameter $\sigma\in(0,1)$ is introduced to make our analysis extendable.
In what follows, we use subscripted $C$, such as $C_{\gamma}$, $C_v$ and $C_u$, to denote a generic positive constant,
which is not necessarily the same at different occurrences, yet is always
dependent on the given data and the solution but independent of temporal and spatial mesh sizes.

\subsection{Discrete fractional Gr\"{o}nwall lemma}\label{sec:gronwall}

We recall the recent developed discrete fractional Gr\"{o}nwall inequality that involves the well-known Mittag--Leffler function $E_\alpha(z):=\sum_{k=0}^\infty\frac{z^k}{\Gamma(1+k\alpha)}$
in \cite[Lemma~2.2, Theorems~3.1-3.2]{Liao2018discrete}.
\begin{lemma}\label{lem:FractGronwall}
For $n=1,2,\cdots,N$, assume that the discrete convolution kernels $\{A_{n-k}^{(n)}\}_{k=1}^n$ satisfy the following two assumptions:\\
\textbf{Ass1}. There is a constant $\pi_{A}>0$ such that $A_{n-k}^{(n)}\geq \frac{1}{\pi_{A}\tau_{k}}\int_{t_{k-1}}^{t_{k}}\omega_{1-\alpha}(t_{n}-s)\zd s$ for $1\leq k\leq n$.\\
\textbf{Ass2}. The discrete kernels are monotone, i.e. $A_{n-k-1}^{(n)}\ge A_{n-k}^{(n)}$ for $1\leq k\leq n-1$.\\
We define a sequence of discrete complementary convolution kernels~$\{P_{n-j}^{(n)}\}_{j=1}^n$ by
\begin{align}\label{eq: discreteConvolutionKernel-RL}
P_{0}^{(n)}:=\frac{1}{A_{0}^{(n)}},\quad
P_{n-j}^{(n)}:=
\frac{1}{p_0^{(j)}}
\sum_{k=j+1}^{n}\brab{A_{k-j-1}^{(k)}-A_{k-j}^{(k)}}P_{n-k}^{(n)},
    \quad 1\leq j\leq n-1.
\end{align}
Then the discrete complementary convolution kernels $P^{(n)}_{n-j}\ge0$ fulfill
\begin{align}
&\sum_{j=k}^nP^{(n)}_{n-j}A_{j-k}^{(j)}\equiv1,
\quad\text{for $1\le k\le n\le N$.}\label{eq: P A}\\
&\sum_{j=1}^nP^{(n)}_{n-j}\omega_{1+m\alpha-\alpha}(t_j)\leq \pi_{A}\omega_{1+m\alpha}(t_n),
\quad\text{for $m=0,1$ and $1\le n\le N$.}\label{eq: P bound}
\end{align}
Suppose that $\lambda_{0}$ and $\lambda_{1}$ are non-negative constants independent of the time-steps, $\lambda:=\lambda_{0}+\lambda_{1}$
and the maximum step size $\tau\le1/\sqrt[\alpha]{2\Gamma(2-\alpha)\lambda\pi_{A}}.$
If the non-negative sequences $(v^k)_{k=0}^N$, $(\xi^{k})_{k=1}^{N}$ and $(\eta^{k})_{k=1}^{N}$ satisfy
\begin{equation}\label{eq: first Gronwall}
\sum_{k=1}^nA_{n-k}^{(n)}\triangledown_{\tau} v^k\le
\lambda_{0}v^{n}
+\lambda_{1}v^{n-1}
+\xi^n+\eta^n\quad\text{for\ $1\le n\le N$,}
\end{equation}
then for $1\le n\le N$ it holds that
\begin{align*}
v^n&\le2E_\alpha\big(2\max\{1,\rho\}\lambda\pi_{A} t_{n}^{\alpha}\big)
	\Big(v^0+\max_{1\le k\le n}\sum_{j=1}^kP^{(k)}_{k-j}(\xi^j+\eta^{j})\Big)\\
&\le2E_\alpha\big(2\max\{1,\rho\}\lambda\pi_{A} t_{n}^{\alpha}\big)
	\Big(v^0+\Gamma(1-\alpha)\pi_{A}\max_{1\le k\le n}\{t_k^{\alpha}\xi^{k}\}+\pi_{A}\omega_{1+\alpha}(t_n)\max_{1\le k\le n}\eta^{k}\Big).
\end{align*}
\end{lemma}

\section{A second-order maximum principle preserving scheme}
In this section, we shall present our second order fully discrete scheme for the time-fractional Allen-Cahn equation (\ref{Problem-1})-(\ref{Problem-2}). In what follows, we consider $\theta:=\alpha/2.$
\subsection{The Alikhanov formula under nonuniform grids}
Given a grid function $\{v^{k}\}$ that is defined on a nonuniform grid (\ref{nonuniformgrid}), for $k\geq{1},$ we define the difference operator $\triangledown_{\tau}v^{k}:=v^{k}-v^{k-1}$,
the difference quotient operator $\partial_{\tau}v^{k-\frac12}:=\triangledown_{\tau}v^{k}/\tau_k$
and the weighted operator $v^{k-\theta}:=(1-\theta)v^{k}+\theta v^{k-1}$.
We then denote by $\Pi_{1,k}v$ the linear interpolant of a function $v$ with respect to the nodes $t_{k-1}$ and $t_{k}$,
and by $\Pi_{2,k}v$ the quadratic with respect to the nodes $t_{k-1}, t_{k}$ and $t_{k+1}$.
The corresponding interpolation errors are denoted by
$\brab{\widetilde{\Pi_{\nu,k}}v}(t):=v(t)-\bra{\Pi_{\nu,k}v}(t)$ for $\nu=1,2$.

Recalling that $\rho_{k}=\tau_{k}/\tau_{k+1}$, then it is easy to show (by using the Newton form of the interpolating polynomials) that
\[
\brab{\Pi_{1,k}v}^{\prime}(t)
=\frac{\triangledown_{\tau}v^{k}}{\tau_{k}}\quad
\text{and}\quad
\brab{\Pi_{2,k}v}^{\prime}(t)
=\frac{\triangledown_{\tau}v^{k}}{\tau_{k}}
+\frac{2(t-t_{k-1/2})}{\tau_{k}(\tau_{k}+\tau_{k+1})}
\braB{\rho_{k}\triangledown_{\tau}v^{k+1}-\triangledown_{\tau}v^{k}}.
\]
The nonuniform Alikhanov approximation \cite{Liao2016JSC,Liao2018second} to~$(\partial_{t}^{\alpha}v)(t_{n-\theta})$ is defined by
\begin{align}\label{eq: L2-1_sigma}
(\partial_{\tau}^{\alpha}v)^{n-\theta}
	:=&\,\int_{t_{n-1}}^{t_{n-\theta}}\omega_{1-\alpha}(t_{n-\theta}-s)\bra{\Pi_{1,n}v}'(s)\zd{s}
		+\sum_{k=1}^{n-1}\int_{t_{k-1}}^{t_k}\omega_{1-\alpha}(t_{n-\theta}-s)\bra{\Pi_{2,k}v}'(s)\zd{s}\nonumber\\
   =&\,a^{(n)}_0\diff v^n+\sum_{k=1}^{n-1}\braB{
	a^{(n)}_{n-k}\diff v^k+\rho_k b^{(n)}_{n-k}\diff v^{k+1}
		-b^{(n)}_{n-k}\diff v^k}.
\end{align}
Here and hereafter, we set  $\sum_{k=i}^{j}\cdot = 0 $ if $i>j$.
The associated discrete convolution kernels $a_{n-k}^{(n)}$ and $b_{n-k}^{(n)}$
are defined, respectively, as
\begin{align}
a_{n-k}^{(n)}
&:=\frac{1}{\tau_{k}}\int_{t_{k-1}}^{\min\{t_{k},t_{n-\theta}\}}
\omega_{1-\alpha}(t_{n-\theta}-s)\zd{s},\quad 1\leq{k}\leq{n};\label{convolution-Kernels: a}\\
b_{n-k}^{(n)}
&:=\frac{2}{\tau_{k}(\tau_{k}+\tau_{k+1})}
\int_{t_{k-1}}^{t_{k}}(s-t_{k-\frac{1}{2}})\omega_{1-\alpha}(t_{n-\theta}-s)\zd{s},\quad
1\leq{k}\leq{n-1}.\label{convolution-Kernels: b}
\end{align}
By re-organizing the terms in~\eqref{eq: L2-1_sigma} we obtain the following compact form
\begin{align}\label{Alikhanov-Formula}
\bra{\partial_{\tau}^{\alpha}v}^{n-\theta}
:=\sum_{k=1}^{n}A_{n-k}^{(n)}\triangledown_{\tau}v^{k},
\end{align}
where the discrete kernels $A_{n-k}^{(n)}$ are defined by:
$A_{0}^{(1)}:=a_{0}^{(1)}$ if $n=1$ and for $n\geq{2}$,
\begin{align}\label{Convolution-Kernels}
A_{n-k}^{(n)}:=\begin{cases}
a_{0}^{(n)}+\rho_{n-1}b_{1}^{(n)},&\text{for}\;k=n,\vspace{2mm}\\
a_{n-k}^{(n)}+\rho_{k-1}b_{n-k+1}^{(n)}-b_{n-k}^{(n)},
&\text{for}\; 2\leq{k}\leq{n-1},\vspace{2mm}\\
a_{n-1}^{(n)}-b_{n-1}^{(n)},&\text{for}\; k=1.
\end{cases}
\end{align}
Notice that the above nonuniform formula is an extension of the Alikhanov Formula on the uniform mesh \cite{Alikhanov2015A}, where the positiveness and monotonicity of $A_{n-k}^{(n)}=A_{n-k}$ were established.
The nonuniform version (\ref{Alikhanov-Formula}) was first proposed in \cite{Liao2016JSC} to resolve the initial singularity by using
a graded mesh near the initial time. Recently, the following results are presented in \cite[Theorem 2.2]{Liao2018second}:

\begin{lemma}\label{lem:Coefficient-Estimate}
Let \emph{\textbf{M1}} hold and consider the discrete convolution kernels $A_{n-k}^{(n)}$ in \eqref{Convolution-Kernels}.
\begin{itemize}
  \item [(i)] The discrete kernels $A_{n-k}^{(n)}$ fulfill $A^{(n)}_{0}\leq\frac{24}{11\tau_{n}}\int_{t_{n-1}}^{t_{n}}\omega_{1-\alpha}(t_{n}-s)\zd{s}$
  and
  \[
  A^{(n)}_{n-k}\geq\frac{4}{11\tau_{n}}\int_{t_{n-1}}^{t_{n}}\omega_{1-\alpha}(t_{n}-s)\zd{s},
  \quad 1\leq{k}\leq{n};\]
  \item [(ii)] The discrete kernels $A_{n-k}^{(n)}$ are monotone for $1\leq{k}\leq{n-1}$,
     \[
     A^{(n)}_{n-k-1}-A^{(n)}_{n-k}
     \geq (1+\rho_{k})b^{(n)}_{n-k}-\frac{1}{5\tau_{k}}
     \int_{t_{k-1}}^{t_{k}}(t_{k}-s)\omega_{-\alpha}(t_{n-\theta}-s)\zd{s}>0.
     \]
  \item [(iii)] And the first kernel $A^{(n)}_{0}$ is appropriately larger than the second one,
      \[
      \frac{1-2\theta}{1-\theta}A^{(n)}_{0}-A^{(n)}_{1}>0 \quad for \quad n\geq{2}.
      \]
\end{itemize}
\end{lemma}
We remark that the estimates in Lemma \ref{lem:Coefficient-Estimate} are much more stronger than
the previous results in \cite{Alikhanov2015A,Liao2016JSC} on the uniform mesh, and these estimates will play an important role when analyzing our adaptive time stepping schemes for phase field equations (e.g., the Allen-Cahn equation in this work). In particular, the boundedness and monotonicity of $A_{n-k}^{(n)}$ are essential to verify the discrete maximum principle
of our second-order time-stepping scheme for the Allen-Cahn equation.

Lemma \ref{lem:Coefficient-Estimate} also implies that the discrete convolution kernels $A_{n-k}^{(n)}$
satisfy the two assumptions \textbf{Ass1}-\textbf{Ass2} in Lemma \ref{lem:FractGronwall} with $\pi_{A}=\frac{11}{4}$, and this will be adapted to show the convergence analysis of our time-stepping scheme using  the discrete complementary convolution kernel argument.

\subsection{The second order fully discrete scheme}
To present the fully discrete scheme, we recall briefly the  difference approximation in physical domain.
For a positive integer $M$,
let the spatial length $h:=L/M$.
Also, we denote $\bar{\Omega}_{h}:=\big\{\mathbf{x}_{h}=(ih,jh)\,|\,0\leq i,j\leq M\}$
and set $\Omega_{h}:=\bar{\Omega}_{h}\cap\Omega$.
For any grid function $\{v_h\,|\,\mathbf{x}_{h}\in\bar{\Omega}_{h}\}$, we denote the grid function space as
\[
\mathbb{V}_{h}:=\big\{v\,|\,v=(v_{j})^{T}\;\;\text{for}\;\;1\leq{j}\leq{M},
\;\;\text{with}\;\;
v_{j}=(v_{i,j})^{T}\;\;\text{for}\;\;1\leq{i}\leq{M}\big\},
\]
where $v^{T}$ is the transpose of the vector $v$.
The maximum norm $\|v\|_{\infty}$ is defined as $\|v\|_{\infty}:=\max_{\mathbf{x}_{h}\in\Omega_{h}}|v_{h}|$.

We shall use the center difference scheme for discetizing the Laplace operator $\Delta$ subject to periodic boundary conditions.
To this end, we denote by $D_{h}$ the associated discrete matrix, then we have $D_{h}=I\otimes{D}+D\otimes{I}$ with $\otimes$ being the Kronecker tensor product operator and
\[
D
=\frac{1}{h^{2}}
\left(
\begin{array}{ccccc}
-2 &  1 & 0 & \cdots & 1 \\
1  & -2 & 1 & \cdots & 0 \\
\vdots & \ddots & \ddots & \ddots & \vdots\\
0  & \cdots & 1 & -2 & 1 \\
1  & \cdots & 0 & 1  &-2 \\
\end{array}
\right)_{M\times{M}}.
\]
We are now ready to present our time-weighted difference scheme for {\eqref{Problem-1}}-{\eqref{Problem-2}}:
\begin{align}
\brab{\partial_{\tau}^{\alpha}u}^{n-\theta}
&=\varepsilon^{2}D_{h}u^{n-\theta}
-f(u)^{n-\theta},\quad
{n}\geq{1},\label{Scheme-1}\\
u_{h}^{0}&=u_{0}(\mathbf{x}_{h}),\quad
\mathbf{x}_{h}\in\bar{\Omega}_{h},\label{Scheme-2}
\end{align}
where the weighted nonlinear term $f(u)^{n-\theta}$ is given by
$$f(u)^{n-\theta}:=\theta{f}(u^{n-1})+(1-\theta)f(u^{n}),$$
and the vector $f(u^{n})$ is defined in the element-wise: $f(u^{n}):=(u^{n})^{.3}-u^{n}$.

To show the uniquely solvability of the above scheme, we list some useful properties of the matrix $D_{h}:$
\begin{lemma}\label{lem:Negative-Condition}
The discrete matrix $D_{h}$ has the following properties
\begin{itemize}
  \item [(a)] The discrete matrix $D_{h}$ is symmetric.
  \item [(b)] For any nonzero $v\in{\mathbb{V}_{h}}$, $v^{T}D_{h}v\leq{0}$, i.e., the matrix $D_{h}$ is negative semi-definite.
  \item [(c)] The elements of $D_{h}=(d_{ij})$ fulfill $d_{ii}=-\max_{i}\sum_{j\neq{i}}|d_{ij}|$ for each $i$.
\end{itemize}
\end{lemma}
The above properties are standard results and are easy to verify. We are now ready to show the following lemma.
\begin{lemma}\label{lem:solvability}
The nonlinear difference scheme \eqref{Scheme-1}-\eqref{Scheme-2} is uniquely solvable
if the step-ratio restriction \textbf{M1} holds with the maximum step size $\tau\le \sqrt[\alpha]{\frac{\omega_{2-\alpha}(1-\theta)}{(1-\theta)}}$.
\end{lemma}
\begin{proof}We rewrite the nonlinear scheme \eqref{Scheme-1} into
\begin{align*}
G_hu^{n}+(1-\theta)(u^{n})^{.3}=g(u^{n-1}),\quad
{n}\geq{1},
\end{align*}
where $G_h:=A_0^{(n)}-1+\theta-(1-\theta)\varepsilon^{2}D_{h}$ and
\begin{align*}
g(u^{n-1}):=&\,
\theta\varepsilon^{2}D_{h}u^{n-1}-\theta{f}(u^{n-1})+
\sum_{k=1}^{n-1}\brab{A_{n-k-1}^{(n)}-A_{n-k}^{(n)}}u^k+A_{n-1}^{(n)}u^0
,\quad{n}\geq{1}.
\end{align*}
If $\tau\le \sqrt[\alpha]{\frac{\omega_{2-\alpha}(1-\theta)}{(1-\theta)}}$, then by the definitions \eqref{Convolution-Kernels}
and \eqref{convolution-Kernels: a} we have
 \begin{align}\label{lower bound A0}
A_0^{(n)}\ge a_0^{(n)}
=\frac{\omega_{2-\alpha}(1-\theta)}{\tau_n^{\alpha}}\ge 1-\theta.
\end{align}
Thus the matrix $G_h$ is positive definite according to Lemma \ref{lem:Negative-Condition} (b).
Consequently, the solution of the nonlinear equations solves
  \begin{align*}
u^n=\arg \min_{w\in \mathbb{V}_{h}}\left\{\frac12w^TG_hw+\frac{1-\theta}4\sum_{k=1}^{M}w_k^{4}-w^Tg(u^{n-1})\right\}\quad
\text{for $n\geq1$.}
\end{align*}
The strict convexity of the above objective function implies the unique solvability of \eqref{Scheme-1}-\eqref{Scheme-2}. The proof is completed.
\end{proof}
\subsection{Discrete maximum principle}
In this section, we show the discrete maximum principle for our scheme \eqref{Scheme-1}-\eqref{Scheme-2}. To this end, we first recall the following lemma \cite[Lemma3.2]{Hou2017Numerical}.
\begin{lemma}\label{lem:Matrix-Inf-Norm}
Let $B$ be a real $M\times{M}$ matrix and $A=aI-B$ with $a>0$.
If the elements of $B=(b_{ij})$ fulfill $b_{ii}=-\max_{i}\sum_{j\neq{i}}|b_{ij}|$,
then for any $c>0$ and $V\in{\mathbb{R}^{M}}$ we have
\begin{align}
\|AV\|_{\infty}\geq{a}\|V\|_{\infty}\quad\text{and}\quad
\|AV+c(V)^{3}\|_{\infty}\geq{a}\|V\|_{\infty}+c\|V\|_{\infty}^{3}.\nonumber
\end{align}
\end{lemma}
Now we are ready to establish the following theorem.
using Lemmas \ref{lem:Coefficient-Estimate}, \ref{lem:Negative-Condition} and \ref{lem:Matrix-Inf-Norm}.
\begin{theorem}\label{thm:Dis-Max-Principle}
Assume that the ratio restriction \textbf{M1} holds and the maximum step size
\begin{align}\label{Time-Step-Constraint}
\tau\leq\min
\bigg\{\sqrt[\alpha]{\frac{\theta\omega_{2-\alpha}(1-\theta)}{2(1-\theta)}},
\sqrt[\alpha]{\frac{h^2\omega_{2-\alpha}(1-\theta)}{4\varepsilon^{2}}}\bigg\}\,.
\end{align}
The second-order scheme {\eqref{Scheme-1}}-{\eqref{Scheme-2}}
preserves the maximum principle \eqref{Frac-Maximum-Principle} at the discrete levels and is unconditionally stable,
that is,  for $1\leq{k}\leq{N}$ we have
$$\mynormb{u^{k}}_{\infty}\leq{1}, \quad \textmd{if} \quad  \mynormb{u^{0}}_{\infty}\leq{1}. $$
\end{theorem}
\begin{proof}
We shall use the mathematical induction argument.
Obviously, the claimed inequality holds for $n=0$.
For $1\leq n\le N$, assume that
\begin{align}\label{inductionAssume}
\mynormb{u^{k}}_{\infty}\leq{1}\quad\text{for $0\leq{k}\leq{n-1}.$}
\end{align}
It remains to verify that $\mynormb{u^{n}}_{\infty}\leq{1}$.
From the definition \eqref{Alikhanov-Formula}, we have
\begin{align*}
(\partial_{\tau}^{\alpha}u)^{n-\theta}
=A_{0}^{(n)}u^{n}
-\brab{A_{0}^{(n)}-A_{1}^{(n)}}u^{n-1}-\mathcal{L}^{n-2}(u),
\end{align*}
where $\mathcal{L}^{n-2}(u)$ is given by
\begin{align}\label{histroy operator}
\mathcal{L}^{n-2}(u):=\sum_{k=1}^{n-2}\big(A_{n-k-1}^{(n)}-A_{n-k}^{(n)}\big)u^{k}+A_{n-1}^{(n)}u^{0}.
\end{align}
Then the scheme \eqref{Scheme-1} can be formulated as follows
\begin{align}\label{Scheme-Variant}
(A^{(n)}_{0}-1+\theta)u^{n}
-&\,(1-\theta)\varepsilon^{2}D_{h}u^{n}
+(1-\theta)(u^{n})^{.3}\nonumber\\
=&\,\brab{A_{0}^{(n)}-A_{1}^{(n)}}u^{n-1}+\theta\varepsilon^{2}D_{h}u^{n-1}
+\theta\kbra{u^{n-1}-(u^{n-1})^{.3}}+\mathcal{L}^{n-2}(u)\nonumber\\
=&\,Q_hu^{n-1}+\theta\kbra{\brab{A_{0}^{(n)}-A_{1}^{(n)}+1}u^{n-1}-(u^{n-1})^{.3}}+\mathcal{L}^{n-2}(u),
\end{align}
where the matrix $Q_h$ is defined by
\begin{align}\label{matrix Qh}
Q_h:=(1-\theta)\brab{A_{0}^{(n)}-A_{1}^{(n)}}+\theta\varepsilon^{2}D_{h}.
\end{align}

We first handle the first term at the right hand side of \eqref{Scheme-Variant}.
It is easy to check that the matrix $Q_h=(q_{ij})$ satisfies $q_{ij}\ge0$ for $i\neq j$, and
\begin{align*}
q_{ii}=(1-\theta)\brab{A_{0}^{(n)}-A_{1}^{(n)}}-\frac{4\theta\varepsilon^{2}}{h^2}
\quad\text{and}\quad
\max_{i}\sum_{j}q_{ij}\le (1-\theta)\brab{A_{0}^{(n)}-A_{1}^{(n)}}.
\end{align*}
Assuming that $\tau\le \sqrt[\alpha]{\frac{h^2}{4\varepsilon^{2}}\omega_{2-\alpha}(1-\theta)}$,
then by Lemma \ref{lem:Coefficient-Estimate} (iii) and \eqref{lower bound A0} we obtain
\begin{align*}
(1-\theta)\brab{A_{0}^{(n)}-A_{1}^{(n)}}
>\theta A^{(n)}_{0}\ge
\frac{\theta}{\tau_n^{\alpha}}\omega_{2-\alpha}(1-\theta)\ge\frac{4\theta\varepsilon^{2}}{h^2},
\end{align*}
or $q_{ii}\ge0$.  Thus all elements of $Q_h$ are nonnegative and
\begin{align*}
\mynormb{Q_h}_{\infty}=\max_{i}\sum_{j}\abs{q_{ij}}=\max_{i}\sum_{j}q_{ij}\le (1-\theta)\brab{A_{0}^{(n)}-A_{1}^{(n)}}.
\end{align*}
Consequently, the induction hypothesis \eqref{inductionAssume} yields
\begin{align}\label{estimate1-Scheme-Variant}
\mynormb{Q_hu^{n-1}}_{\infty}\le \mynormb{Q_h}_{\infty}\mynormb{u^{n-1}}_{\infty}\le (1-\theta)\brab{A_{0}^{(n)}-A_{1}^{(n)}}.
\end{align}

For the second term at the right hand side of of \eqref{Scheme-Variant}, consider the following function
\begin{align*}
\psi(z)
:=\brab{A_{0}^{(n)}-A_{1}^{(n)}+1}z-z^{3}.
\end{align*}
If $\tau\le \sqrt[\alpha]{\frac{\theta}{2(1-\theta)}\omega_{2-\alpha}(1-\theta)}$,
Lemma \ref{lem:Coefficient-Estimate} (iii) and \eqref{lower bound A0} give
\begin{align*}
A_{0}^{(n)}-A_{1}^{(n)}
>\frac{\theta}{1-\theta} A^{(n)}_{0}\ge \frac{\theta\omega_{2-\alpha}(1-\theta)}{(1-\theta)\tau_n^{\alpha}}\ge2.
\end{align*}
In this case, one has $\abs{\psi(z)}\le A_{0}^{(n)}-A_{1}^{(n)}$ for any $z\in[-1,1]$. Therefore,
the induction hypothesis \eqref{inductionAssume} yields
\begin{align}\label{estimate2-Scheme-Variant}
\theta\mynormb{\brab{A_{0}^{(n)}-A_{1}^{(n)}+1}u^{n-1}-(u^{n-1})^{.3}}_{\infty}\le \theta\brab{A_{0}^{(n)}-A_{1}^{(n)}}.
\end{align}
For the last term $\mathcal{L}^{n-2}(u)$ of \eqref{Scheme-Variant}, the decreasing property in Lemma \ref{lem:Coefficient-Estimate} (ii)
 and the induction hypothesis \eqref{inductionAssume} lead to
\begin{align}\label{estimate3-Scheme-Variant}
\mynormb{\mathcal{L}^{n-2}(u)}_{\infty}\leq \sum_{k=1}^{n-2}\brab{A_{n-k-1}^{(n)}-A_{n-k}^{(n)}}\mynormb{u^{k}}_{\infty}
+A_{n-1}^{(n)}\mynormb{u^{0}}_{\infty}\leq A_{1}^{(n)}.
\end{align}
Moreover, under the setting $\tau\le \sqrt[\alpha]{\frac{\omega_{2-\alpha}(1-\theta)}{1-\theta}}$,
the inequality \eqref{lower bound A0} shows $A_0^{(n)}>1-\theta$.
Then by using Lemmas \ref{lem:Negative-Condition} and \ref{lem:Matrix-Inf-Norm}, one can
bound the left hand side of \eqref{Scheme-Variant} by
\begin{align*}
&\quad\mynormb{(A^{(n)}_{0}-1+\theta)u^{n}-(1-\theta)\varepsilon^{2}D_{h}u^{n}+(1-\theta)(u^{n})^{.3}}_{\infty}\\
&\ge (A^{(n)}_{0}-1+\theta)\mynormb{u^{n}}_{\infty}+(1-\theta)\mynormb{u^{n}}_{\infty}^3.
\end{align*}
Therefore, collecting the estimates \eqref{estimate1-Scheme-Variant}--\eqref{estimate3-Scheme-Variant},
it follows from \eqref{Scheme-Variant} that
\begin{align*}
&\qquad (A^{(n)}_{0}-1+\theta)\,\mynormb{u^{n}}_{\infty}+(1-\theta)\mynormb{u^{n}}_{\infty}^3\\
&\le\,\mynorm{Q_hu^{n-1}+\theta
\kbrab{\brab{A_{0}^{(n)}-A_{1}^{(n)}+1}u^{n-1}-(u^{n-1})^{.3}}+\mathcal{L}^{n-2}(u)}_{\infty}\\
&\le\,\mynormb{Q_hu^{n-1}}_{\infty}
+\theta\mynormb{\brab{A_{0}^{(n)}-A_{1}^{(n)}+1}u^{n-1}-(u^{n-1})^{.3}}_{\infty}
+\mynormb{\mathcal{L}^{n-2}(u)}_{\infty}\\
&\le\, (1-\theta)\brab{A_{0}^{(n)}-A_{1}^{(n)}}+\theta\brab{A_{0}^{(n)}-A_{1}^{(n)}}+A_{1}^{(n)}=A^{(n)}_{0}.
\end{align*}
This immediately implies $\mynormb{u^{n}}_{\infty}\leq1$.
Otherwise, we have
$$(A^{(n)}_{0}-1+\theta)\mynormb{u^{n}}_{\infty}
+(1-\theta)\mynormb{u^{n}}_{\infty}^3-A_{0}^{(n)}>0,$$
as the function
$
\phi(z):=(A^{(n)}_{0}-1+\theta)z
+(1-\theta)z^{3}-A_{0}^{(n)}
$
is monotonically increasing for any $z>0$. This leads to a contradiction and the proof is completed.
\end{proof}

We remark that the maximum time-step restriction \eqref{Time-Step-Constraint} is only a sufficient condition to
ensure the discrete maximum principle (see Example \ref{example:maximum principle}).
In the time-fractional Allen-Cahn equation \eqref{Problem-1},
the coefficient $\varepsilon\ll 1$ represents the width of diffusive interface. Always,
we should choose a small space length $h=O(\varepsilon)$ to track the moving interface.
So, in most situations, the restriction \eqref{Time-Step-Constraint} is practically reasonable because
it is approximately equivalent to
$$\tau\leq\sqrt[\alpha]{\frac{\theta\omega_{2-\alpha}(1-\theta)}{2(1-\theta)}}\rightarrow \frac12\quad\text{as $\alpha\rightarrow1$.}$$
Notice also that the condition \eqref{Time-Step-Constraint} may become worse when the fractional order $\alpha\rightarrow0$.
However, this time-step condition is sharp in the sense that
it is compatible with the restriction in \cite{Hou2017Numerical} that ensures the discrete maximum principle of
Crank-Nicolson scheme for the integer-order Allen-Cahn equation.

\section{Error convolution structure and convergence analysis}

We consider the error analysis by denoting the consistency error of Alikhanov formula \eqref{Alikhanov-Formula} as
$\Upsilon^{j}[v]:=(\partial_{t}^{\alpha}v)(t_{j-\theta})-(\partial_{\tau}^{\alpha}v)^{j-\theta}$ for $j\ge1$.
Similar as in \cite[Theorem 3.4]{Liao2018second}, we show in the next lemma that $\Upsilon^{j-\theta}$ can be controlled by a discrete convolution structure,
which is valid for a general class of time meshes.
Moreover, the fractional Gr\"{o}nwall inequality in Lemma \ref{lem:FractGronwall}
suggests that the solution error is determined by the convolution error $\sum_{j=1}^{n}P_{n-j}^{(n)}\absb{\Upsilon^{j}[v]}$,
where $P_{n-j}^{(n)}$ are the discrete complementary convolution kernels
defined in \eqref{eq: discreteConvolutionKernel-RL}.

\begin{lemma}\label{lem:Alikhanov-Formula-Consistence}
Assume that the step-ratio condition \emph{\textbf{M1}} holds,
the function $v\in{C^{3}((0,T])}$ and $\int_0^Ts^2\abs{v'''(s)}\zd{s}<\infty$.
For the nonuniform Alikhanov formula \eqref{Alikhanov-Formula} with
the discrete convolution kernels $A_{n-k}^{(n)}$, the local consistency error $\Upsilon^{j-\theta}$
has a convolution structure
\begin{align*}
\absb{\Upsilon^{n}[v]}
\leq
A_{0}^{(n)}G_{\emph{loc}}^{n}
+\sum_{k=1}^{n-1}\brab{A_{n-k-1}^{(n)}-A_{n-k}^{(n)}}G_{\emph{his}}^{k},
\quad  1\leq{n}\leq{N},
\end{align*}
where the terms $G_{\emph{loc}}^{k}$ and $G_{\emph{his}}^{k}$ are defined by, respectively,
\begin{align*}
G_{\emph{loc}}^{k}&:=
\frac{3}{2}\int_{t_{k-1}}^{t_{k-1/2}}\bra{s-t_{k-1}}^{2}|v'''(s)|\zd{s}
+\frac{3\tau_{k}}{2}\int_{t_{k-1/2}}^{t_{k-1}}\bra{t_{k}-s}|v'''(s)|\zd{s}
\\
G_{\emph{his}}^{k}&:=
\frac{5}{2}\int_{t_{k-1}}^{t_{k}}\bra{s-t_{k-1}}^{2}|v'''(s)|\zd{s}
+\frac{5}{2}\int_{t_{k}}^{t_{k+1}}\bra{t_{k+1}-s}^{2}|v'''(s)|\zd{s}.
\end{align*}
Consequently, the global convolution error satisfies
\begin{align*}
\sum_{j=1}^nP_{n-j}^{(n)}\absb{\Upsilon^{j}[v]}
\leq\sum_{k=1}^{n}P_{n-k}^{(n)}A_{0}^{(k)}G_{\emph{loc}}^{k}
+\sum_{k=1}^{n-1}P_{n-k}^{(n)}A_{0}^{(k)}G_{\emph{his}}^{k},
\quad  1\leq{n}\leq{N}.
\end{align*}
\end{lemma}

Notice that the global consistency error in
Lemma \ref{lem:Alikhanov-Formula-Consistence}
gives a super\/convergence estimate of nonuniform Alikhanov formula.
Consider the first time level $n=1$, the regularity setting
\eqref{Regularity-Sigma} gives
$$\abs{\Upsilon^{1}}\leq
A_{0}^{(1)}G_{\mathrm{loc}}^{1}
\leq{C}_{u}\tau_1^{\sigma-\alpha}/\sigma,$$
which implies $\Upsilon^{1}=O(1)$ when $\sigma=\alpha$,
and if $0<\sigma\le\alpha$ then the situation becomes worse.
However, we have the global consistency error of order $\mathcal{O}(\tau_{1}^{\sigma})$ (see Tables \ref{Error-Test-1}-\ref{Error-Test-2} in Section 5)
as one has $P_{0}^{(1)}\abs{\Upsilon^{1}}
\leq G_{\mathrm{loc}}^{1}\leq{C}_{u}\tau_{1}^{\sigma}/\sigma.$
In general, Lemma \ref{lem:Alikhanov-Formula-Consistence} leads to the following corollary (see also \cite[Lemma 3.6]{Liao2018second}).
\begin{corollary}\label{Global-Consis-Error}
Assume that the step-ratio condition \emph{\textbf{M1}} holds,
and the function $v\in{C^{3}((0,T])}$
admits an initial singularity, $\abs{v'''(t)}\le C_v(1+t^{\sigma-2})$ as $t\rightarrow0$
for a real parameter $0<\sigma<1$.
The global consistency error can be bounded by
\begin{align*}
\sum_{j=1}^{n}P_{n-j}^{(n)}\absb{\Upsilon^{j}[v]}
\leq C_v\braB{\,\tau_1^{\sigma}/\sigma
+t_{1}^{\sigma-3}\tau_{2}^{3}+
+\frac{1}{1-\alpha}\max_{2\leq{k}\leq{n}}
t_{k}^{\alpha}t_{k-1}^{\sigma-3}\tau_{k}^{3}/\tau_{k-1}^{\alpha}}, \quad 1\leq{n}\leq{N}.
\end{align*}
Specifically, if the mesh satisfies the graded-like condition \emph{\textbf{M2}}, then
\begin{align*}
\sum_{j=1}^{n}P_{n-j}^{(n)}\absb{\Upsilon^{j}[v]}
\leq
\frac{C_v}{\sigma(1-\alpha)}\tau^{\min\{\gamma\sigma,2\}},
\quad  1\leq{n}\leq{N}.
\end{align*}
\end{corollary}
The next lemma \cite[Lemma 3.8]{Liao2018second} shows that the temporal error introduced by the time
weighted approximation is bounded by the error that is generated by the Alikhanov approximation.
\begin{lemma}\label{Weighted-Approach-Error}
Assume that $v\in{C^{2}((0,T])}$, and there exists a positive constant $C_{v}$ such that
$\abs{v^{\prime\prime}(t)}\leq{C}_{v}\bra{1+t^{\sigma-2}}$ for $0<t\leq{T}$, where $\sigma\in(0,1)$ is a regularity parameter.
Denote the local truncation error of $v^{n-\theta}$ by
\[
R^{n}[v]=v(t_{n-\theta})-v^{n-\theta},\quad 1\leq{n}\leq{N}.
\]
If the graded-like condition \emph{\textbf{M2}} holds, then the global consistency error satisfies
\[
\sum_{j=1}^{n}P_{n-j}^{(n)}\absb{R^{n}[v]}
\leq{C}_{v}\bra{\tau_{1}^{\sigma+\alpha}/\sigma
+t_{n}^{\alpha}\max_{2\leq{k}\leq{n}}t_{k-1}^{\sigma-2}\tau_{k}^{2}},\quad 1\leq{n}\leq{N}.
\]
\end{lemma}

Taking the advantage of the discrete maximum principle
in Theorem \ref{thm:Dis-Max-Principle}, one can prove the convergence of numerical solution
without assuming the Lipschitz continuity of the nonlinear term $f(u)$.
More precisely, we have the following error estimates.
\begin{theorem}\label{Convergence-Theorem}
Assume that $\mynormb{u^{0}}_{L^{\infty}}\leq{1}$ and
the solution of \eqref{Problem-1}-\eqref{Problem-2} satisfies the regular assumption \eqref{Regularity-Sigma}.
If the ratio restriction \textbf{M1} holds and the maximum step size
\begin{align*}
\tau\leq\min
\bigg\{\sqrt[\alpha]{\frac{\omega_{2-\alpha}(1)}{11}},\,
\sqrt[\alpha]{\frac{\theta\omega_{2-\alpha}(1-\theta)}{2(1-\theta)}},\,
\sqrt[\alpha]{\frac{h^2\omega_{2-\alpha}(1-\theta)}{4\varepsilon^{2}}}\bigg\}\,,
\end{align*}
then the solution of \eqref{Scheme-1}-\eqref{Scheme-2} is convergent in the maximum norm, that is,
\begin{align*}
\mynormb{u(t_n)-u^{n}}_{\infty}
\leq C_u\braB{\frac{\tau_1^{\sigma}}{\sigma}
+\frac{1}{1-\alpha}\max_{2\leq{k}\leq{n}}t_{k}^{\alpha}t_{k-1}^{\sigma-3}\tau_{k}^{3-\alpha}
+h^{2}}, \quad 1\leq{n}\leq{N}.
\end{align*}
Specially, when the time mesh satisfies \emph{\textbf{M2}}, it holds that
\begin{align*}
\mynormb{u(t_n)-u^{n}}_{\infty}
\leq \frac{C_u}{\sigma(1-\alpha)}
\tau^{\min\{\gamma\sigma,2\}}
+C_{u}h^{2},\quad 1\leq{n}\leq{N}.
\end{align*}
Notice that the proposed scheme achieves the optimal accuracy $O(\tau^{2})$ if the graded parameter $\gamma\geq\max{\{1,\,2/\sigma\}}$.
\end{theorem}
\begin{proof} We set $U_{h}^{n}:=u(\mathbf{x}_h,t_n)$ and denote the error function as $e_{h}^{n}:=U_{h}^{n}-u_{h}^{n}\in{\mathbb{V}_{h}}$
for $\mathbf{x}_{h}\in\bar{\Omega}_{h}$ and $0\leq{n}\leq{N}$.
It is easy to find that the exact solution $U_{h}^{n}$ satisfies the governing equations
\begin{align*}
\bra{\partial_{\tau}^{\alpha}U}^{n-\theta}
-\varepsilon^{2}D_{h}U^{n-\theta}
&=-f(U)^{n-\theta}
+\Upsilon^{n}[u]+R^{n}[u]+R_{s}^{n},\quad {1}\leq{n}\leq{N},\\
U_{h}^{0}
&=u_{0}(\mathbf{x}_{h}),\quad \mathbf{x}_{h}\in\Omega_{h},
\end{align*}
where $R_{s}^{n}$ represents the truncation errors in space.
It is easy to get the error equation
\begin{align}
\bra{\partial_{\tau}^{\alpha}e}^{n-\theta}
-\varepsilon^{2}D_{h}e^{n-\theta}
&=-f(U)^{n-\theta}
+f(u)^{n-\theta}+\Upsilon^{n}[u]+R^{n}[u]+R_{s}^{n},\quad {1}\leq{n}\leq{N},\label{Error-Equation-1}
\end{align}
subject to the zero-valued initial data $e^{0}=0$. To facilitate the subsequent analysis,
we rewrite the equation \eqref{Error-Equation-1} into the following form
\begin{align}
A_{0}^{(n)}e^{n}
+\mathcal{L}^{n-2}(e)
-(1-\theta)\varepsilon^{2}D_{h}e^{n}
&= Q_he^{n-1}+\theta\brab{A_{0}^{(n)}-A_{1}^{(n)}}e^{n-1}\nonumber\\
&\quad +f(u)^{n-\theta}
-f(U)^{n-\theta}
+\Upsilon^{n}[u]+R^{n}[u]+R_{s}^{n},\label{Error-Equation-3}
\end{align}
where $\mathcal{L}^{n-2}(e)$  and $Q_h$ are defined by \eqref{histroy operator} and \eqref{matrix Qh}, respectively.
Recalling the inequality
$$\absb{(a^{3}-a)-(b^{3}-b)}\leq{2}\absb{a-b}\quad \text{for $\forall\, a,b\in[-1,1]$},$$
we apply Theorem \ref{thm:Dis-Max-Principle} to obtain
\begin{align*}
\mynormb{f(u)^{n-\theta}-f(U)^{n-\theta}}_{\infty}
\leq
2\theta\mynormb{e^{n-1}}_{\infty}
+2(1-\theta)\mynormb{e^{n}}_{\infty}.
\end{align*}
With the help of the triangle inequality and the estimate \eqref{estimate1-Scheme-Variant}, it follows from \eqref{Error-Equation-3} that
\begin{align}\label{Error-Equation-Estimate-1}
&\quad\,\, \mynormb{A_{0}^{(n)}e^{n}
+\mathcal{L}^{n-2}(e)-(1-\theta)\varepsilon^{2}D_{h}e^{n}}_{\infty}\nonumber\\
&\leq
\brab{A_{0}^{(n)}-A_{1}^{(n)}}\mynormb{e^{n-1}}_{\infty}
+2(1-\theta)\mynormb{e^{n-1}}_{\infty}\nonumber\\
&\quad +2\theta\mynormb{e^{n}}_{\infty}+\mynormb{\Upsilon^{n}[u]}_{\infty}
+\mynormb{R^{n}[u]}_{\infty}
+\mynormb{R_{s}^{n}}_{\infty}.
\end{align}
By using Lemma \ref{lem:Coefficient-Estimate} (iii) and the triangle inequality,
we bound the left hand side of \eqref{Error-Equation-Estimate-1} by
\begin{align}\label{Error-Equation-Estimate-2}
&\quad\,\,\mynormb{A_{0}^{(n)}e^{n}
+\mathcal{L}^{n-2}(e)-(1-\theta)\varepsilon^{2}D_{h}e^{n}}_{\infty}\nonumber\\
&=\mynormB{(A_{0}^{(n)}-\varepsilon^{2}D_{h})e^{n}
-\sum_{k=1}^{n-2}\brab{A_{n-k-1}^{(n)}-A_{n-k}^{(n)}}e^{k}
-A_{0}^{(n)}e^{0}}_{\infty}\nonumber\\
&\geq
A_{0}^{(n)}\mynormb{e^{n}}_{\infty}
-\sum_{k=1}^{n-2}\brab{A_{n-k-1}^{(n)}-A_{n-k}^{(n)}}\mynormb{e^{k}}_{\infty}
-A_{n-1}^{(n)}\mynormb{e^{0}}_{\infty},
\end{align}
where Lemma {\ref{lem:Matrix-Inf-Norm}} and Lemma \ref{lem:Negative-Condition} (c) were used in the last inequality.
Then it follows from \eqref{Error-Equation-Estimate-1}-\eqref{Error-Equation-Estimate-2} that
\begin{align*}
\sum_{k=1}^{n}A_{n-k}^{(n)}\triangledown_{\tau}\mynormb{e^{k}}_{\infty}
\leq2\theta\mynormb{e^{n}}_{\infty}
+2(1-\theta)\mynormb{e^{n-1}}_{\infty}+
\mynormb{\Upsilon^{n}[u]}_{\infty}
+\mynormb{R^{n}[u]}_{\infty}
+\mynormb{R_{s}^{n}}_{\infty},
\end{align*}
which takes the form of \eqref{eq: first Gronwall} with the substitutions $v^k:=\mynormb{e^{k}}_{\infty}$ and
\begin{align*}
\lambda_{0}:=2\theta,\quad
\lambda_{1}:=2(1-\theta),\quad
\xi^{n}:=\mynormb{\Upsilon^{n}[u]}_{\infty}
+\mynormb{R^{n}[u]}_{\infty},\quad
\eta^n:=\mynormb{R_{s}^{n}}_{\infty}.
\end{align*}
Recall that the ratio restriction \textbf{M1} gives $\rho=7/4$ and
Lemma \ref{lem:Coefficient-Estimate} (i) gives $\pi_{A}=\frac{11}{4}$.
The discrete fractional Gr\"{o}nwall inequality in Lemma \ref{lem:FractGronwall} says that, if the maximum time-step size $\tau\leq\sqrt[\alpha]{\frac{\omega_{2-\alpha}(1)}{11}}$, then it holds that
\begin{align*}
\mynormb{e^{n}}_{\infty}\leq
2E_{\alpha}\brab{20t_{n}^{\alpha}}
\bigg[\max_{1\leq{k}\leq{n}}\sum_{j=1}^{k}P_{k-j}^{(k)}\bra{\mynormb{\Upsilon^{j}[u]}_{\infty}+\mynormb{R^{n}[u]}_{\infty}}
+3\omega_{1+\alpha}(t_{n})h^2\bigg].
\end{align*}
Then the desired estimate follows by using together Corollary \ref{Global-Consis-Error} and Lemma \ref{Weighted-Approach-Error}.
\end{proof}

\section{Numerical implementations}
In this section, we provide some details for the numerical implementations.
\subsection{Fast Alikhanov formula}
It is evident that the approximations \eqref{Alikhanov-Formula} is prohibitively expensive for long time simulations due to the long-time memory. Therefore, to reduce the computational cost and storage requirements,
we apply the sum-of-exponentials (SOE) technique to speed up the evaluation of the Alikhanov formula \eqref{Alikhanov-Formula}.
A core result is to approximate the kernel function $\omega_{1-\alpha}(t)$ efficiently on the interval $[\Delta{t},\,T]$, and we shall adopt the results in \cite[Theorem 2.5]{Jiang2017Fast}.
\begin{lemma}\label{SOE}
For the given $\alpha\in(0,\,1)$, an absolute tolerance error $\epsilon\ll{1}$, a cut-off time $\Delta{t}>0$ and a finial time $T$, there exists a positive integer $N_{q}$, positive quadrature nodes $s^{\ell}$ and corresponding positive weights $\varpi^{\ell}\,(1\leq{\ell}\leq{N_{q}})$ such that
\begin{align}
\bigg|
\omega_{1-\alpha}(t)
-\sum_{\ell=1}^{N_{q}}\varpi^{\ell}e^{-s^{\ell}t}\bigg|\leq\epsilon,
\quad
\forall\,{t}\in[\Delta{t},\,T].\nonumber
\end{align}
\end{lemma}
Motivated by the above lemma, we split the Caputo derivative \eqref{CaputoDef} into the sum of a history part (an integral over $[0,\,t_{n-1}]$) and a local part (an integral over $[t_{n-1},\,t_{n}]$) at the time $t_{n}$. Then, the local part will be approximated by linear interpolation directly, the history part can be evaluated via the SOE technique, that is,
\begin{align}\label{SOE-Approximation}
\bra{\partial_{t}^{\alpha}v}\bra{t_{n-\theta}}
&\approx
\int_{t_{n-1}}^{t_{n-\theta}}
\varpi_{n}^{\prime}(s)(\Pi_{1,n}v)'(s)\zd{s}
+\int_{0}^{t_{n-1}}\sum_{\ell=1}^{N_{q}}\varpi^{\ell}
e^{-s^{\ell}(t_{n-\theta}-s)}v^{\prime}(s)\zd{s}\nonumber\\
&=a_{0}^{(n)}\triangledown_{\tau}v^{n}
+\sum_{\ell=1}^{N_{q}}\varpi^{\ell}\mathcal{H}^{\ell}(t_{n-1}),\quad
n\geq{1},
\end{align}
where $\mathcal{H}^{\ell}(t_{0}):=0$ and $\mathcal{H}^{\ell}(t_{k})
:=\int_{0}^{t_{k}}e^{-s^{\ell}(t_{k+1-\theta}-s)}v^{\prime}(s)\zd{s}$.
By using the quadratic interpolation and a recursive formula, we can approximate $\mathcal{H}^{\ell}(t_{k})$ using the following relation
\begin{align}\label{History-Part}
\mathcal{H}^{\ell}(t_{k})
&\approx\int_{0}^{t_{k-1}}e^{-s^{\ell}(t_{k+1-\theta}-s)}v^{\prime}(s)\zd{s}
+\int_{t_{k-1}}^{t_{k}}e^{-s^{\ell}(t_{k+1-\theta}-s)}
(\Pi_{2,k}v)'(s)\zd{s}\nonumber\\
&=e^{-s^{\ell}(\theta\tau_{k}+(1-\theta\tau_{k+1}))}\mathcal{H}^{\ell}(t_{k-1})
+a^{(k,\ell)}\triangledown_{\tau}v^{k}
+b^{(k,\ell)}\brab{\rho_{k}\triangledown_{\tau}v^{k+1}-\triangledown_{\tau}v^{k}},
\end{align}
where the positive coefficients $a^{(k,\ell)}$ and $b^{(k,\ell)}$ are given by, respectively,
\[
a^{(k,l)}:=\frac{1}{\tau_{k}}
\int_{t_{k-1}}^{t_{k}}e^{-s^{\ell}(t_{k+1-\theta}-s)}\zd{s},\quad
b^{(k,l)}:=
\int_{t_{k-1}}^{t_{k}}e^{-s^{\ell}(t_{k+1-\theta}-s)}\frac{2(s-t_{k-1/2})}{\tau_{k}(\tau_{k}+\tau_{k+1})}\zd{s}.
\]
From \eqref{SOE-Approximation}-{\eqref{History-Part}}, we arrive at the fast algorithm of Alikhanov formula
\begin{align}\label{Fast-Approximation}
(\partial_{f}^{\alpha}v)^{n-\theta}
=
a_{0}^{(n)}\triangledown_{\tau}v^{n}
+\sum_{\ell=1}^{N_{q}}\varpi^{\ell}\mathcal{H}^{\ell}(t_{n-1}),
\quad n\geq{1},
\end{align}
in which $\mathcal{H}^{\ell}(t_{k})$ is computed by using the recursive relationship
\begin{align}\label{History-Recursive}
\mathcal{H}^{\ell}(t_{k})
=e^{-s^{\ell}(\theta\tau_{k}+(1-\theta)\tau_{k+1})}\mathcal{H}^{\ell}(t_{k-1})
+a^{(k,\ell)}\triangledown_{\tau}v^{k}
+b^{(k,\ell)}\brab{\rho_{k}\triangledown_{\tau}v^{k+1}-\triangledown_{\tau}v^{k}}.
\end{align}
\subsection{Adaptive time-stepping strategy}
Our theory permitts some adaptive time-stepping strategy to capture the fast dynamics and to reduce the cost of computation.
Roughly speaking, the adaptive time steps can be selected by using an accuracy criterion example as \cite{Gomez2011Provably}, or
the time evolution of the total energy such as \cite{QiaoAn2011}.
We consider the former and update the time step size by using the formula
\begin{align*}
\tau_{ada}\bra{e,\tau}
=S_a\bra{\frac{tol}{e}}^{\frac{1}{2}}\tau,
\end{align*}
where $S_a$ is a default safety coefficient, $tol$ is a reference tolerance, and $e$ is the relative error at each time level.
The adaptive time-stepping strategy is presented in Algorithm \ref{Adaptive-Time-Step-Strategy}.
\begin{algorithm}
\caption{Adaptive time-stepping strategy}
\label{Adaptive-Time-Step-Strategy}
\begin{algorithmic}[1]
\Require{Given $u^{n}$ and time step $\tau_{n}$}
\State Compute $u_{1}^{n+1}$ by a first-order scheme with time step $\tau_{n}$, e.g., the backward Euler-type scheme with $L1$ formula \cite{Liao2018Sharp}.
\State Compute $u_{2}^{n+1}$ by the proposed scheme \eqref{Scheme-1} with time step $\tau_{n}$.
\State Calculate $e_{n+1}=\|u_{2}^{n+1}-u_{1}^{n+1}\|/\|u_{2}^{n+1}\|$.
\If {$e_{n}<tol$ or $\tau_{n}=\frac{2}{3}\tau_{n-1}$}
\State Update time-step size $\tau_{n+1}\leftarrow
\min\{\max\{\tau_{\min},\tau_{ada}\},\tau_{\max}\}$.
\Else
\State Recalculate with time-step size  $\tau_{n}\leftarrow
\max\{\min\{\max\{\tau_{\min},\tau_{ada}\},\tau_{\max}\},\frac{2}{3}\tau_{n-1}\}$.

\State Goto 1
\EndIf
\end{algorithmic}
\end{algorithm}

For the nonlinear time-stepping method \eqref{Scheme-1}-\eqref{Scheme-2}, we adopt an iteration scheme at each time level with the termination error $\eta=10^{-12}$.
The absolute tolerance error of SOE approximation is given as  $\epsilon=10^{-12}$. The maximum norm error $e(N):=\max_{n}\|U^{n}-u^{n}\|_{\infty}$ is recorded in each run,
and the experimental convergence order in time is computed by
\[
\text{Order}\approx\frac{\log\bra{e(N)/e(2N)}}{\log\bra{\tau(N)/\tau(2N)}},
\]
where $\tau(N)$ denotes the maximum time-step size for total $N$ subintervals.

\subsection{Numerical examples}
\begin{table}[htb!]
\begin{center}
\caption{Temporal error of scheme \eqref{Scheme-1}-\eqref{Scheme-2} for $\alpha=0.8,\,\sigma=0.8$ with $\gamma_{\mathrm{opt}}=2.5$.}\label{Error-Test-1} \vspace*{0.3pt}
\def\temptablewidth{1.0\textwidth}
{\rule{\temptablewidth}{0.5pt}}
\begin{tabular*}{\temptablewidth}{@{\extracolsep{\fill}}cccccccccc}
\multirow{2}{*}{$N$} &\multirow{2}{*}{$\tau$} &\multicolumn{2}{c}{$\gamma=1$} &\multirow{2}{*}{$\tau$} &\multicolumn{2}{c}{$\gamma=2.5$} &\multirow{2}{*}{$\tau$}&\multicolumn{2}{c}{$\gamma=4$} \\
             \cline{3-4}          \cline{6-7}         \cline{9-10}
         &          &$e(N)$   &Order &          &$e(N)$   &Order &          &$e(N)$    &Order\\
\midrule
  32     &3.13e-02	 &3.55e-03 &$-$   &7.06e-02	 &4.81e-04 &$-$   &7.95e-02	 &6.80e-04 &$-$\\
  64     &1.56e-02	 &2.04e-03 &0.80  &3.63e-02  &1.19e-04 &2.10  &3.70e-02	 &1.43e-04 &2.04\\
  128    &7.81e-03	 &1.17e-03 &0.80  &1.96e-02  &3.15e-05 &2.15  &2.05e-02	 &3.74e-05 &2.27\\
  256    &3.91e-03	 &6.72e-04 &0.80  &9.20e-03  &5.50e-06 &2.31  &1.04e-02	 &7.68e-06 &2.34\\
\midrule
\multicolumn{3}{l}{$\min\{\gamma\sigma,2\}$}   &0.80 & & &2.00 & & &2.00\\
\end{tabular*}
{\rule{\temptablewidth}{0.5pt}}
\end{center}
\end{table}	
	 	
\begin{table}[htb!]
\begin{center}
\caption{Temporal error of scheme \eqref{Scheme-1}-\eqref{Scheme-2} for $\alpha=0.8,\,\sigma=0.4$ with $\gamma_{\mathrm{opt}}=5$.}\label{Error-Test-2} \vspace*{0.3pt}
\def\temptablewidth{1.0\textwidth}
{\rule{\temptablewidth}{0.5pt}}
\begin{tabular*}{\temptablewidth}{@{\extracolsep{\fill}}cccccccccc}
\multirow{2}{*}{$N$} &\multirow{2}{*}{$\tau$} &\multicolumn{2}{c}{$\gamma=3$} &\multirow{2}{*}{$\tau$} &\multicolumn{2}{c}{$\gamma=5$} &\multirow{2}{*}{$\tau$}&\multicolumn{2}{c}{$\gamma=6$} \\
             \cline{3-4}          \cline{6-7}         \cline{9-10}
         &           &$e(N)$   &Order &          &$e(N)$   &Order &          &$e(N)$   &Order\\
\midrule
  32     &6.85e-02	 &5.87e-03 &$-$   &8.77e-02	 &2.37e-03 &$-$   &8.46e-02	 &2.37e-03 &$-$\\
  64     &3.93e-02	 &2.63e-03 &1.45  &4.32e-02  &6.05e-04 &1.93  &4.56e-02	 &6.07e-04 &2.21\\
  128    &1.91e-02	 &1.16e-03 &1.13  &2.04e-02  &1.51e-04 &1.85  &2.16e-02	 &1.40e-04 &1.96\\
  256    &9.12e-03	 &5.07e-04 &1.12  &1.05e-02  &3.84e-05 &2.08  &1.04e-02	 &3.14e-05 &2.06\\
\midrule
\multicolumn{3}{l}{$\min\{\gamma\sigma,2\}$}   &1.20 & & &2.00 & & &2.00\\
\end{tabular*}
{\rule{\temptablewidth}{0.5pt}}
\end{center}
\end{table}	


\begin{example}\label{example:Accuracy-Test}
We first test the accuracy and consider
$\partial_{t}^{\alpha}u=\varepsilon^{2}\Delta u-f(u)+g(\mathbf{x},t)$
on the space-time domain $(0,1)^{2}\times(0,1]$. We set $\varepsilon=\sqrt{2}/(4\pi)$ and choose an exterior force $g$ such that the exact solution yields $u=\omega_{1+\sigma}(t)\sin(2\pi{x})\sin(2\pi{y})$.
\end{example}
We examine the temporal accuracy using a fine spatial grid mesh with $M=1024$ such that the temporal error dominates the spatial error.
Always, the time interval $[0,T]$ is divided into two parts $[0, T_0]$ and $[T_0, T]$ with total $N$ subintervals.
We will take $T_0=\min\{1/\gamma,T\}$,
and apply the graded grids $t_{k}=T_{0}(k/N_0)^{\gamma}$
in $[0,T_{0}]$ to resolve the initial singularity.
In the remainder interval $[T_{0},T]$,
we put $N_1:=N-N_0$ small cells with random time-step sizes
$\tau_{N_{0}+k}=(T-T_{0})\epsilon_{k}/\sum_{k=1}^{N_1}\epsilon_{k}$ for $1\leq k\leq N_1$,
where $\epsilon_{k}\in(0,1)$ are the random numbers.
The numerical results for two different cases $\sigma=\alpha$ and $\sigma<\alpha$ are listed in Tables \ref{Error-Test-1}-\ref{Error-Test-2}.
It is noticed that the scheme admits a $\mathcal{O}\bra{\tau^{\min\{\gamma\sigma,2\}}}$-order rate of convergence, and thus
the optimal second-order accuracy is achieved when $\gamma\geq\gamma_{opt}=2/\sigma$.

\begin{example}\label{example:Merger-Drops}
We next consider an example of merging of four-drops to show the effectiveness of the adaptive strategy
and to exploit the effect of the fraction order $\alpha$ on the equilibration process. More precisely, we consider
$\partial_{t}^{\alpha}u=\varepsilon^{2}\Delta u-f(u)$ on $\Omega=(-1,1)^{2}\times(0,T]$
with $\varepsilon=0.02$. The solution is computed with $h=0.02$ using the following initial data
\begin{align*}
u_{0}
=&-0.9\tanh\bra{\bra{(x-0.3)^{2}+y^{2}-0.2^2}/\varepsilon}
\tanh\bra{\bra{(x+0.3)^{2}+y^{2}-0.2^2}/\varepsilon}\nonumber\\
&\times\tanh\bra{\bra{x^{2}+(y-0.3)^{2}-0.2^2}/\varepsilon}
\tanh\bra{\bra{x^{2}+(y+0.3)^{2}-0.2^2}/\varepsilon}.
\end{align*}
\end{example}

\begin{figure}[htb!]
\centering
\includegraphics[width=2.7in,height=1.8in]{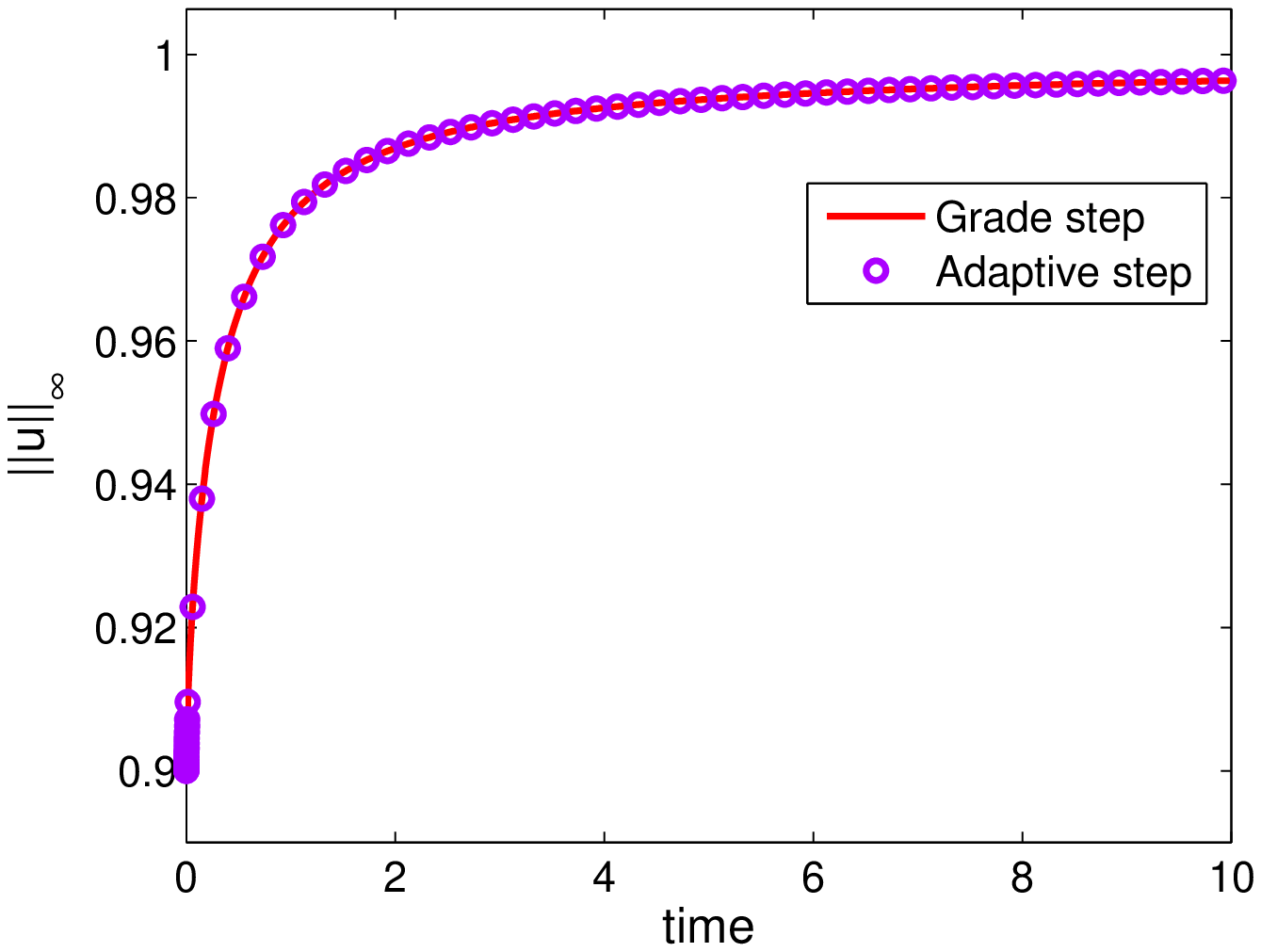}
\includegraphics[width=2.7in,height=1.8in]{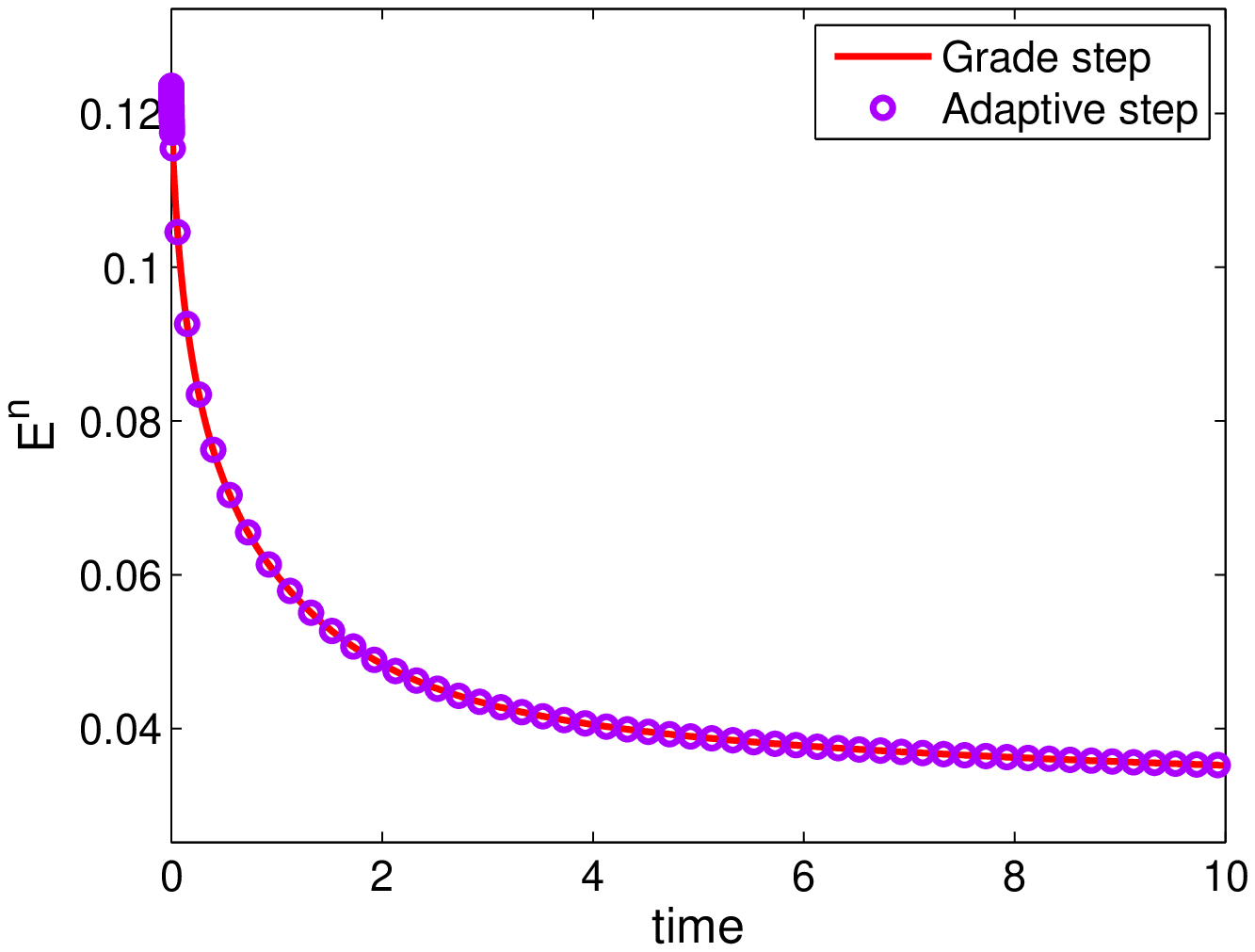}
\caption{The solution in the maximum norm (left) and the discrete energy (right) vary against time until $T=10$
for Example \ref{example:Merger-Drops} with the fractional order $\alpha=0.7$.}
\label{Adaptive-Maximum-Principle-Energy}
\end{figure}

For a fixed fractional order $\alpha=0.7$, Figure \ref{Adaptive-Maximum-Principle-Energy} presents the solution in the maximum norm and
the energy functional $E(t_k)$ over the time interval $[0,T]$ with $T=10.$
The graded mesh with $\gamma=,3$ $N_{0}=30$ and $T_{0}=0.01$ in the starting interval $[0,T_{0}]$ is used to resolve the initial singularity. For $(T_0,T]$ we
first consider a uniform mesh with the total grid number $N_1=970$ (listed as Grade step).
For comparison, we also consider an adaptive grids (listed as Adaptive step), and we use the adaptive time-stepping technique
in the time interval $(T_0,T]$ with the parameters $S_a=0.9$, $tol=10^{-3}$,
and $\tau_{\min}=\tau_{N_{0}}=10^{-3}$ and $\tau_{\max}=10^{-1}$.
It is learned in Figure \ref{Adaptive-Maximum-Principle-Energy} that the adaptive mesh provides good agreement with a fine uniform mesh.
While the adaptive time-stepping strategy leads to a substantial decrease in the computational cost since the number of adaptive steps is 108, while the uniform mesh needs 970 steps.


\begin{figure}[htb!]
\centering
\includegraphics[width=2.7in,height=1.8in]{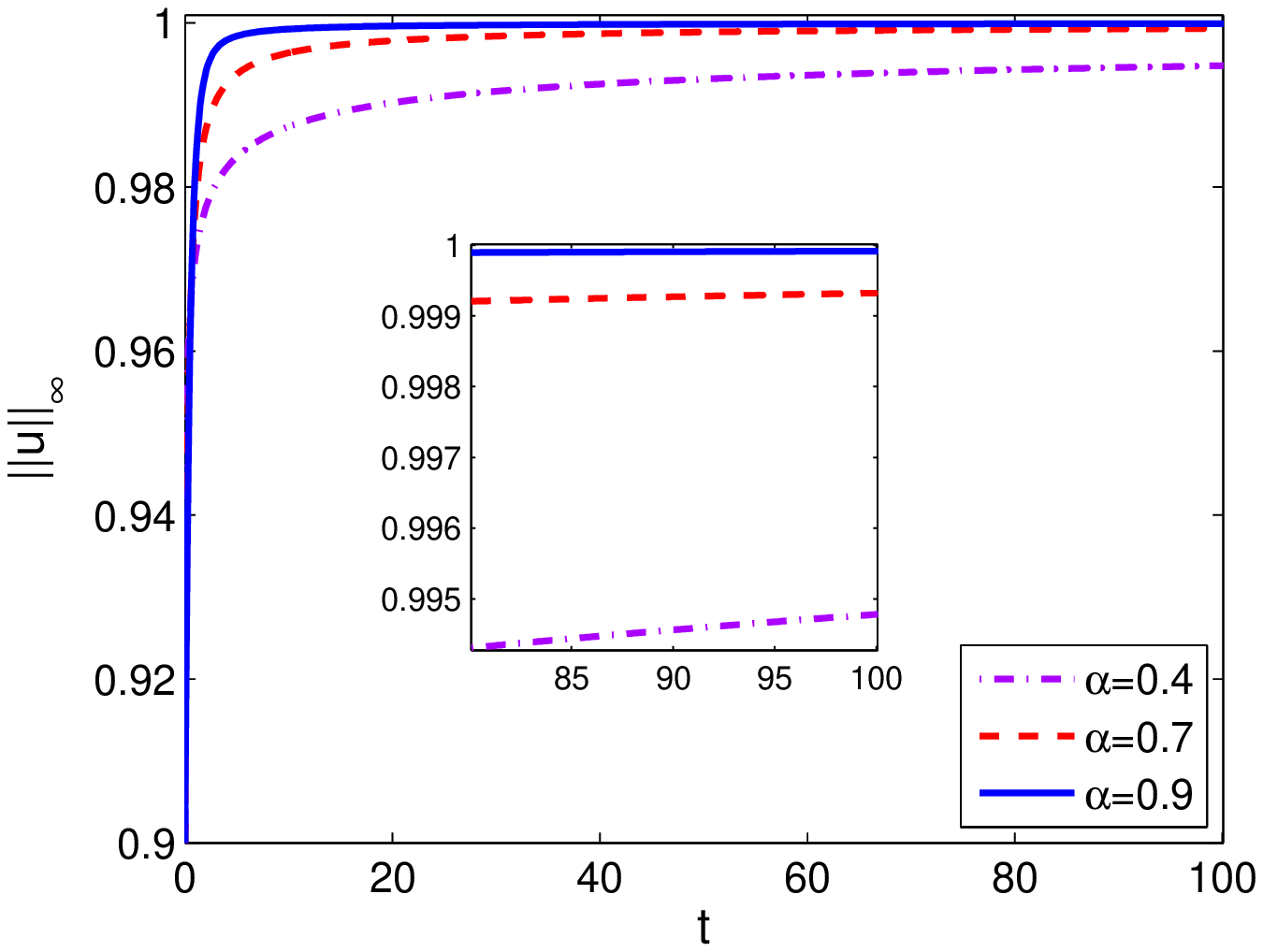}
\includegraphics[width=2.7in,height=1.8in]{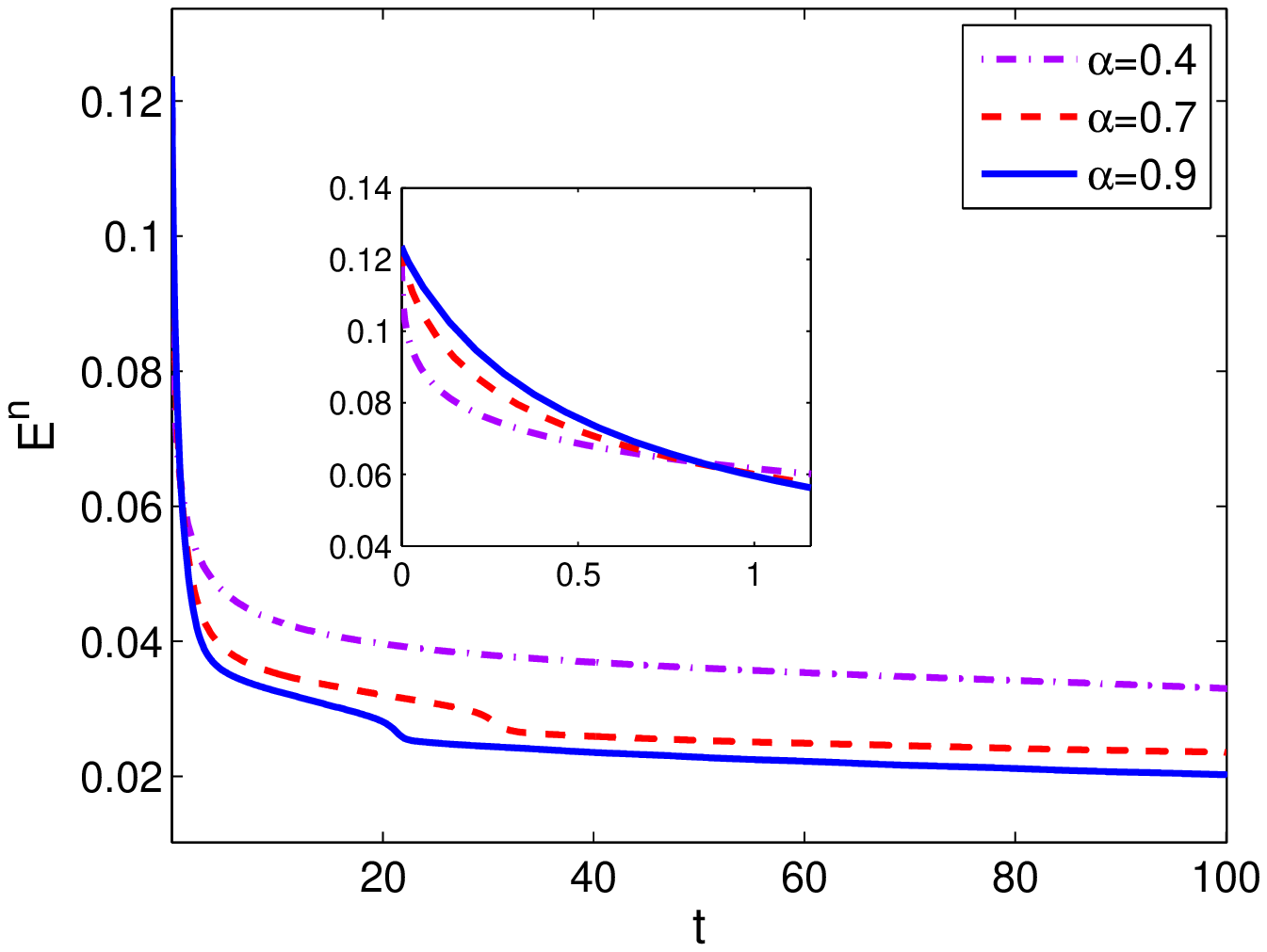}
\caption{The maximum norm values (left) and the discrete energies (right) vary against time until $T=100$
for Example \ref{example:Merger-Drops} with three fractional orders $\alpha=0.4,0.7$ and $0.9$.}
\label{Maximum-Principle-Energy}
\end{figure}

\begin{figure}[htb!]
\centering
\includegraphics[width=1.47in]{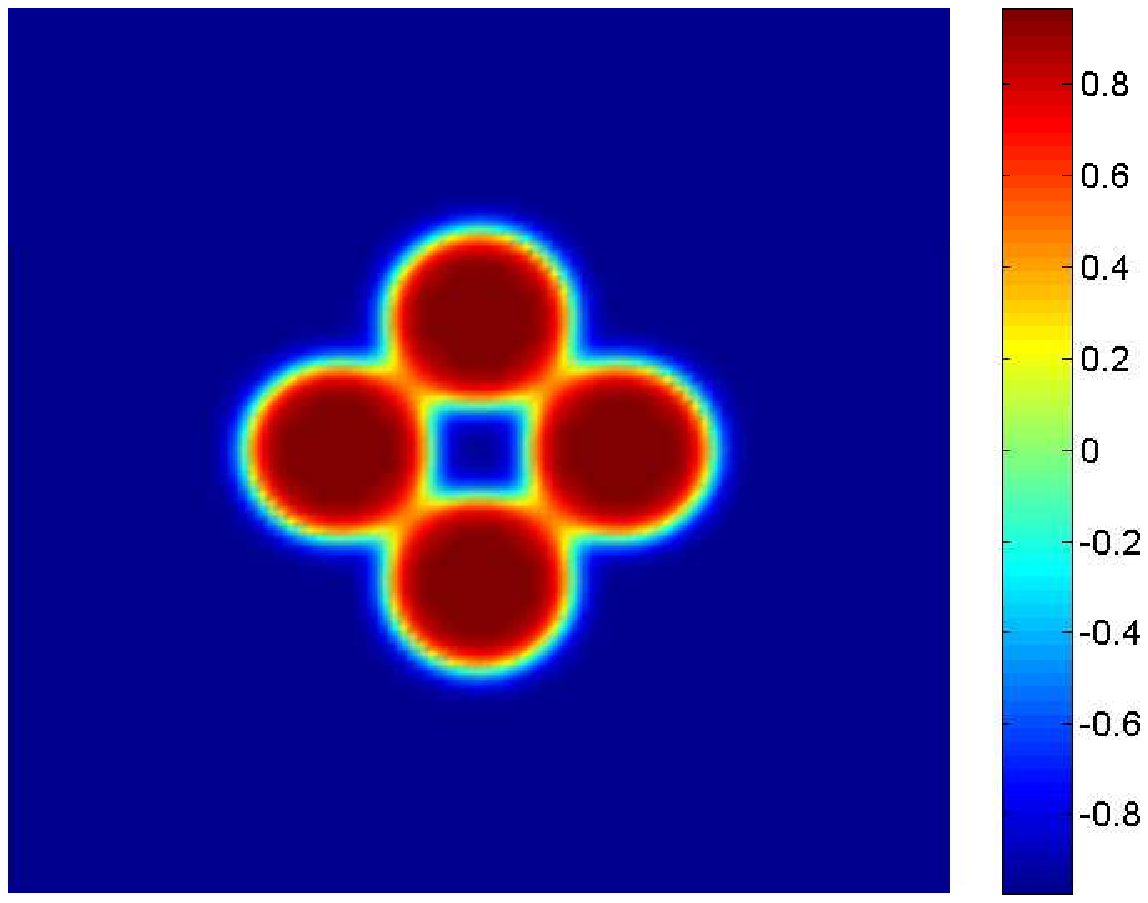}
\includegraphics[width=1.47in]{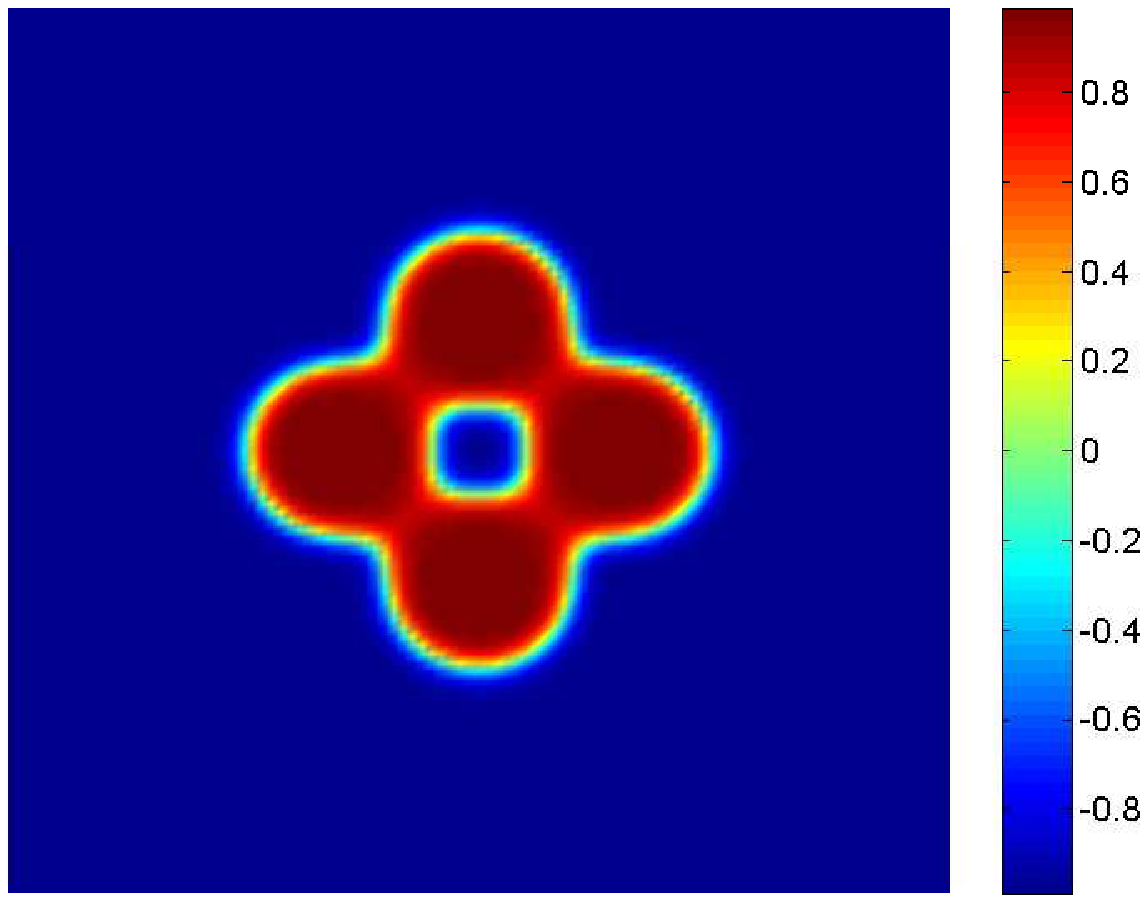}
\includegraphics[width=1.47in]{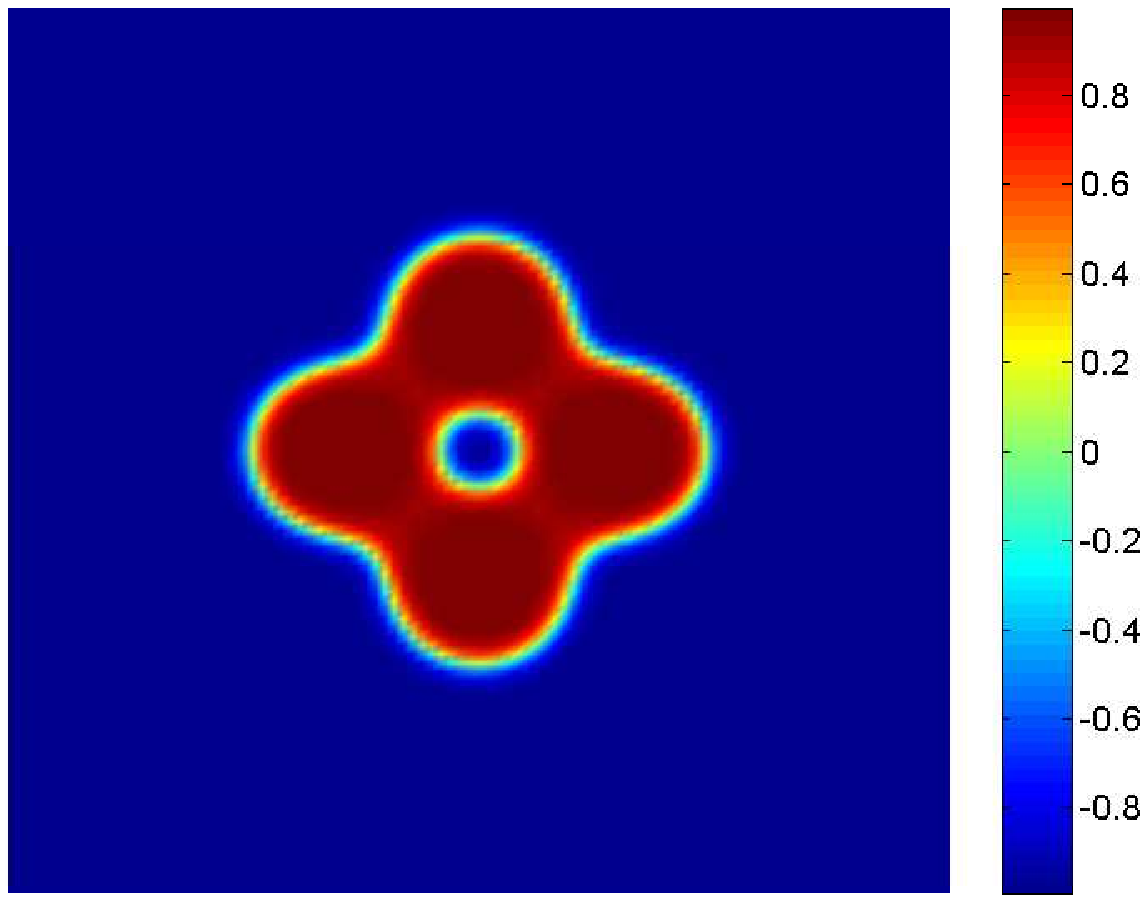}
\includegraphics[width=1.47in]{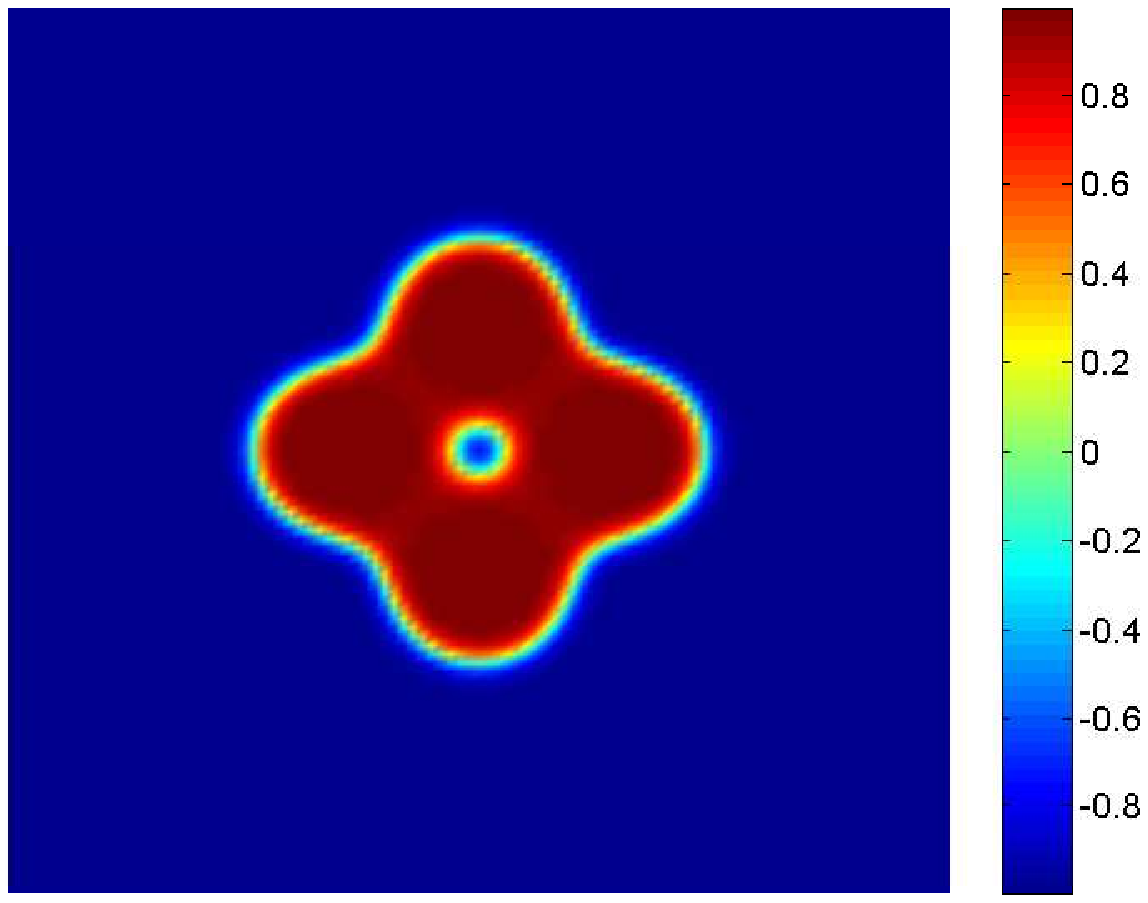}\\
\includegraphics[width=1.47in]{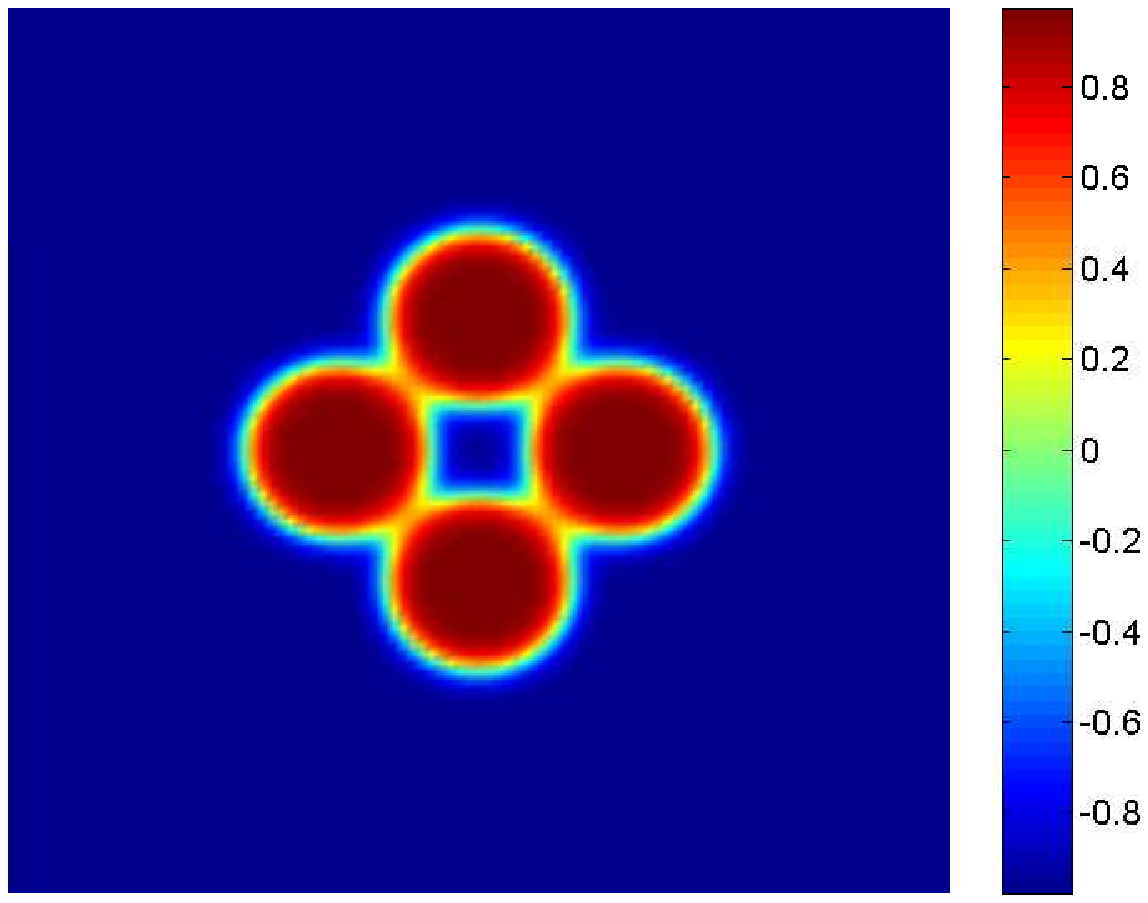}
\includegraphics[width=1.47in]{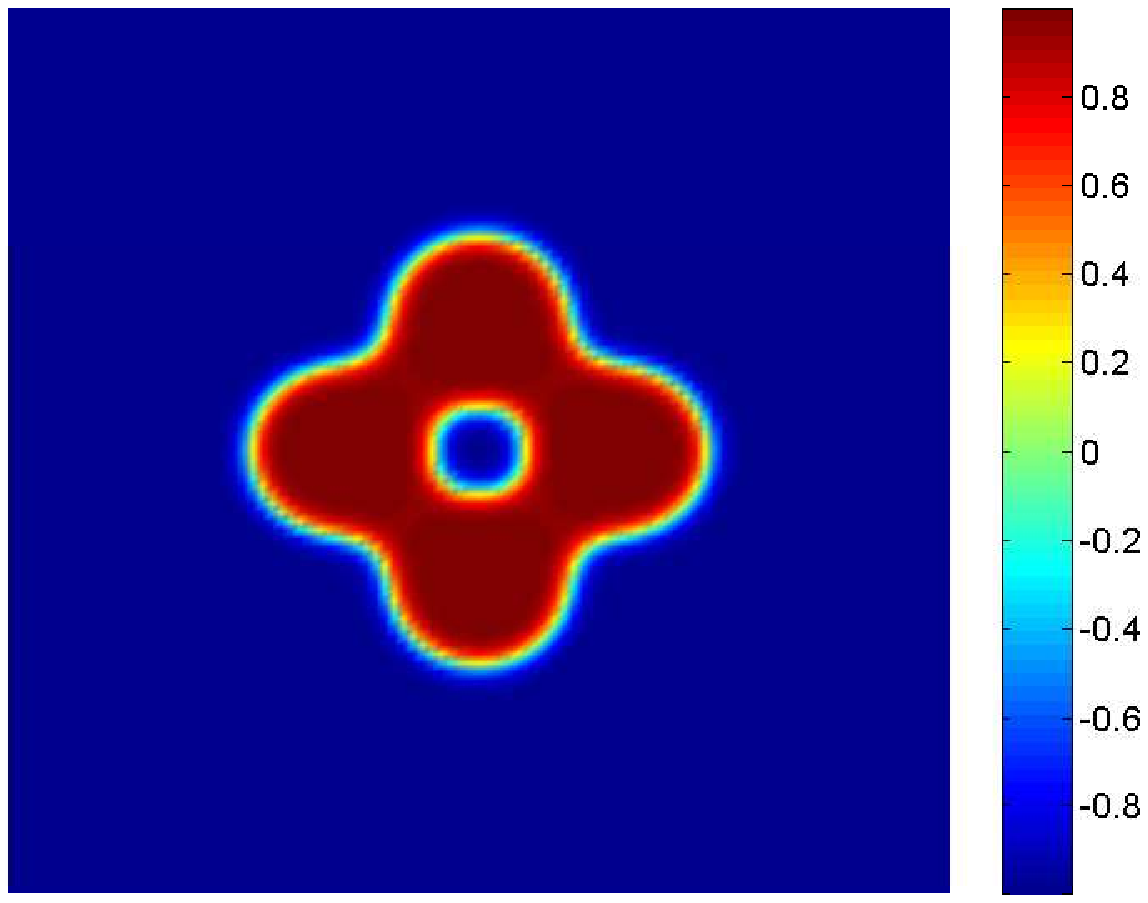}
\includegraphics[width=1.47in]{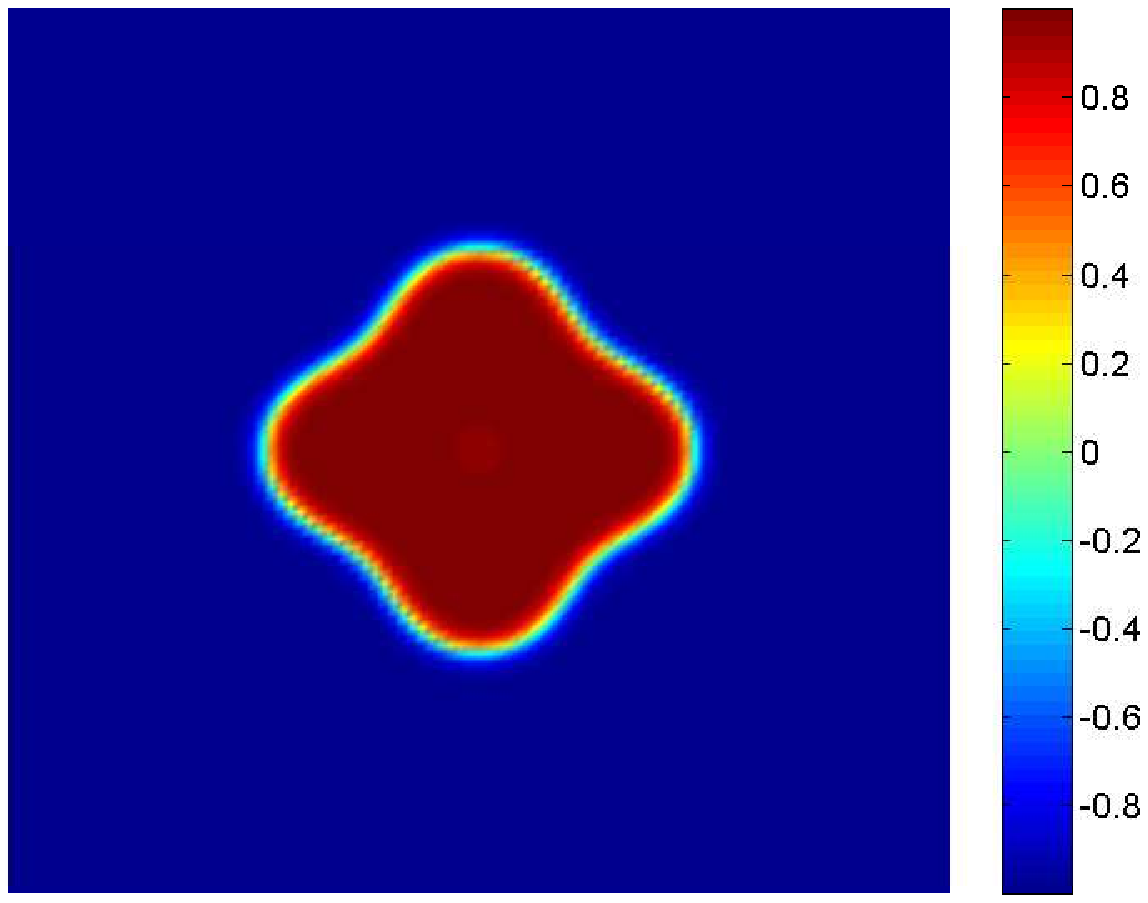}
\includegraphics[width=1.47in]{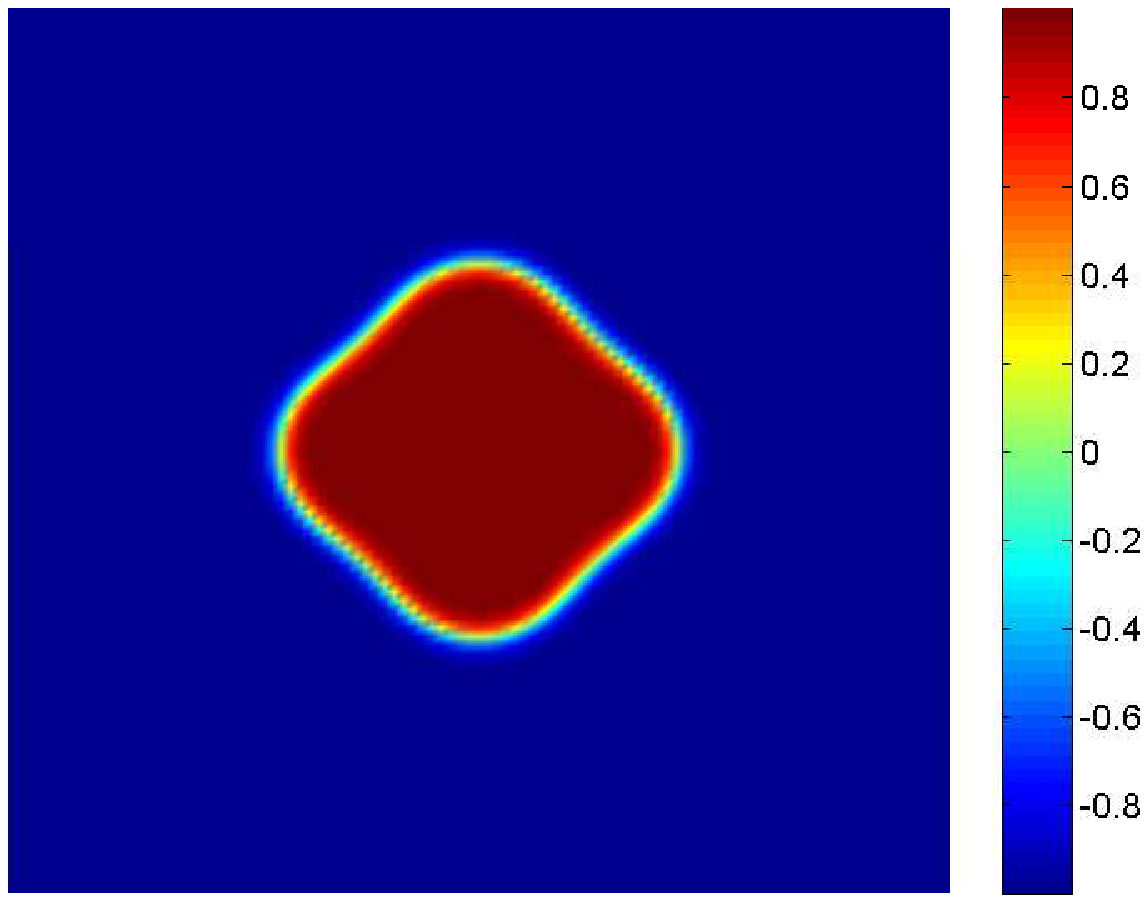}\\
\includegraphics[width=1.47in]{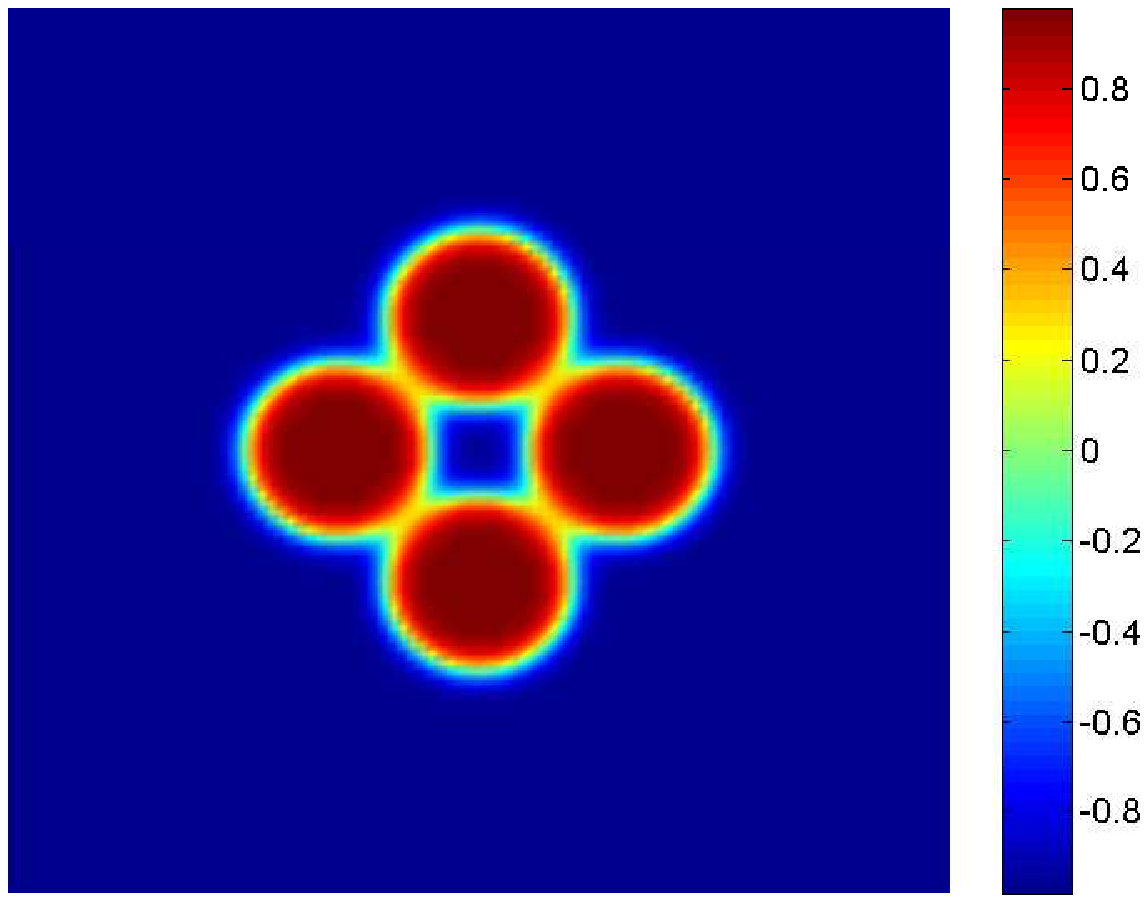}
\includegraphics[width=1.47in]{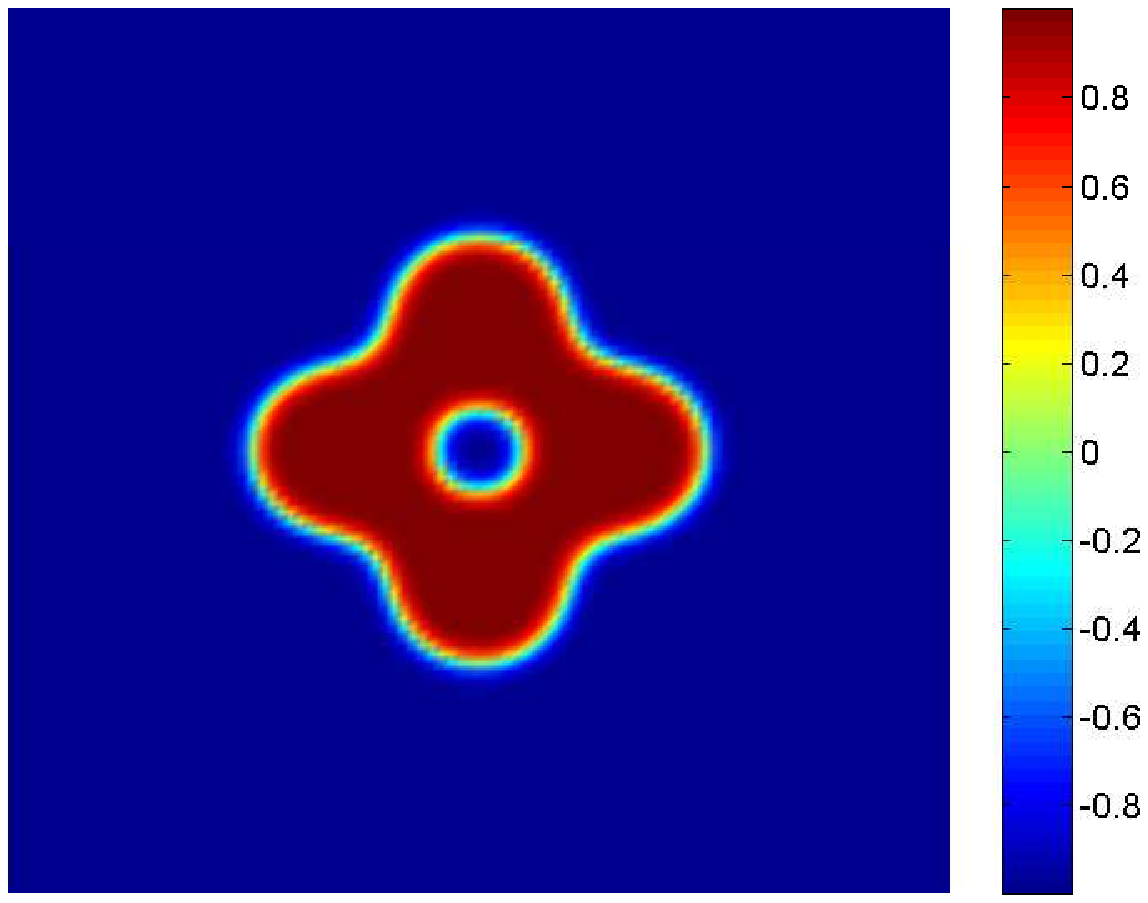}
\includegraphics[width=1.47in]{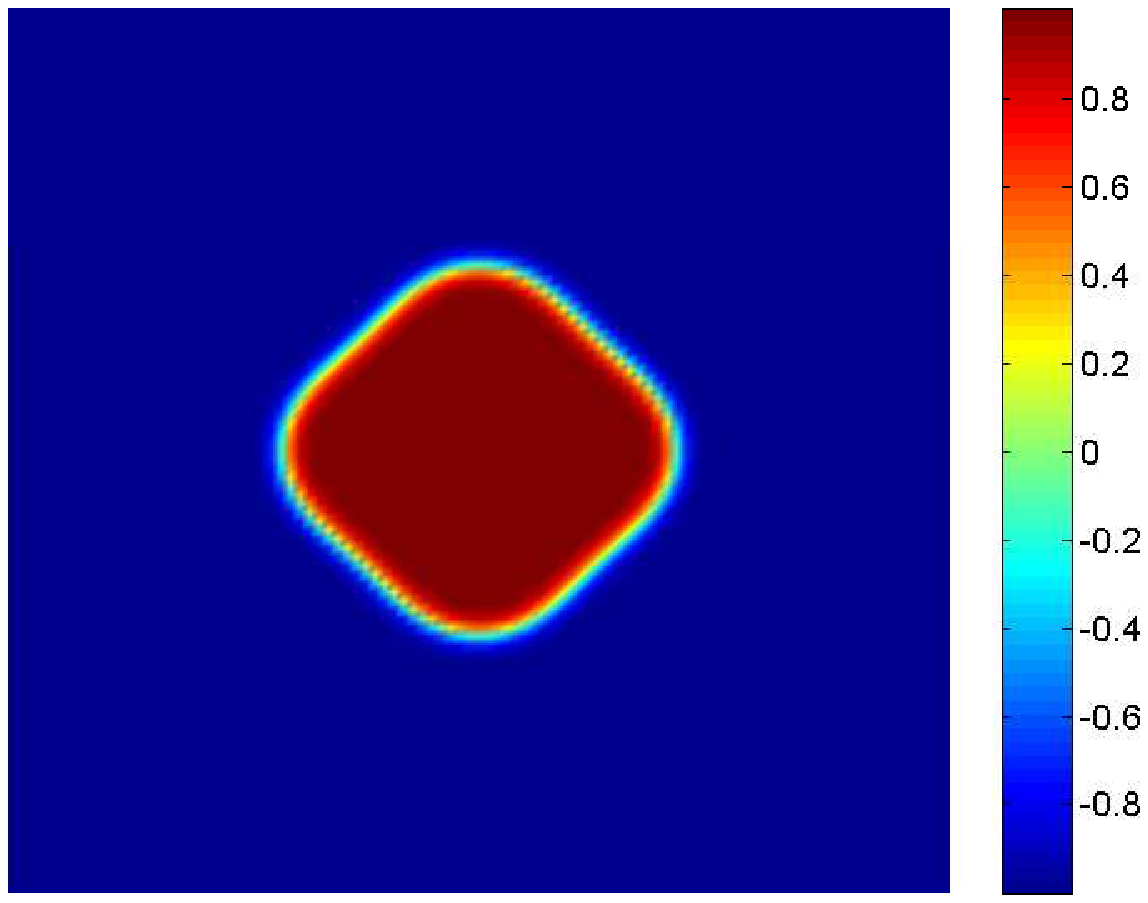}
\includegraphics[width=1.47in]{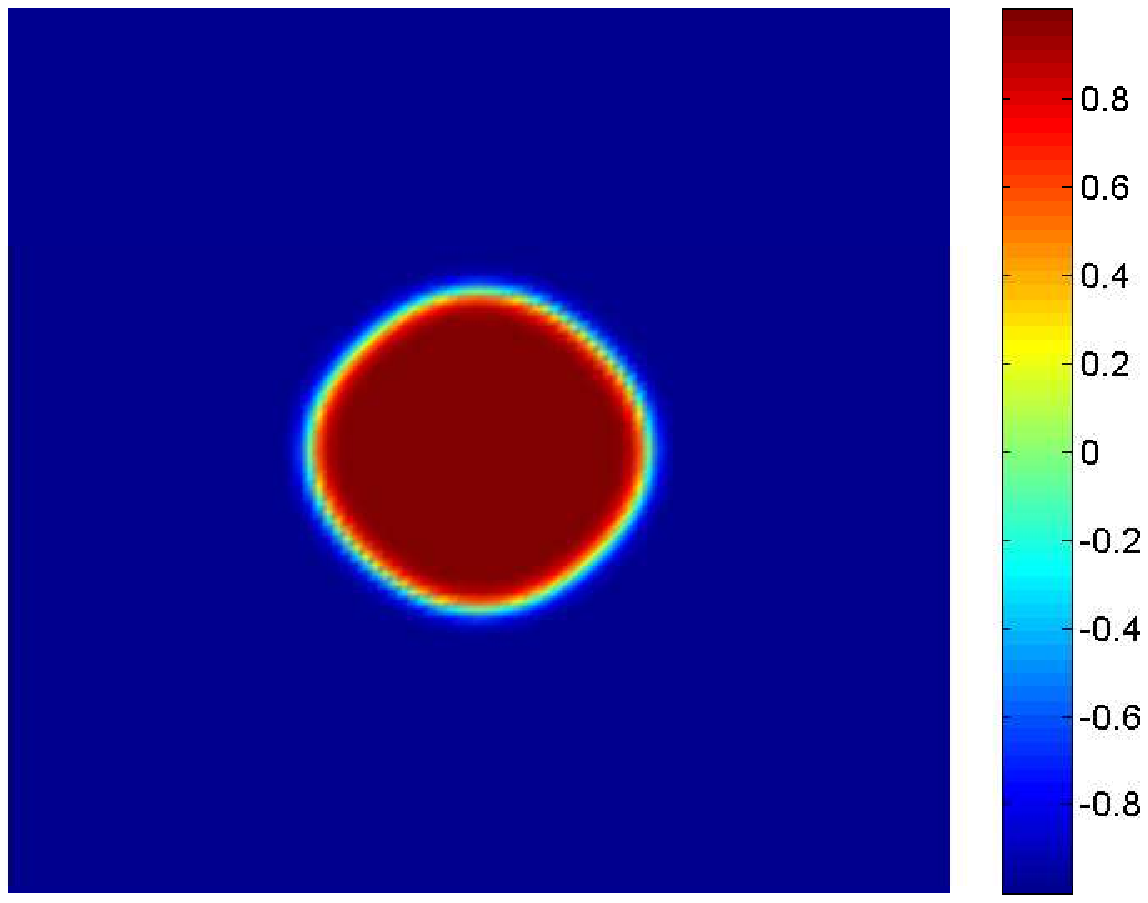}
\caption{Solution snapshots of time-fractional Allen-Cahn equation at $t=1, 10, 50, 100$ (from left to right)
for fractional orders $\alpha=0.4$, $0.7$ and $0.9$ (from top to bottom).}
\label{Merger-Drops-Process}
\end{figure}	 	

Second, we investigate the equilibration process of the drops in Example \ref{example:Merger-Drops}
by using the adaptive strategy. Figure \ref{Maximum-Principle-Energy} compares the maximum norm values and
the discrete energy functionals for three different fractional orders $\alpha=0.4,0.7$ and $0.9$
over a long-time interval $[0,100].$
We observe that the larger the fractional order $\alpha$, the faster the maximum norm value approaches 1,
but the maximum norm values are always bounded by 1 for all cases.  Similarly, the larger the fractional order $\alpha$, the faster the energy dissipates.

Figure \ref{Merger-Drops-Process} displays the snapshots of the solution contours for
different fractional orders $\alpha=0.4,0.7$ and $0.9$. The same adaptive time-stepping technique
is employed in the time interval $(T_0,T]$ with the parameters $S_a=0.9$, $tol=10^{-3}$,
$\tau_{\min}=\tau_{N_{0}}=10^{-3}$ and $\tau_{\max}=10^{-1}$.
As time escapes, the four-drops merges into a single drop
and shrinks progressively (due to the primitive problem dose not conserve the volume).
Moreover, the larger the fractional order $\alpha$, the bigger the shrinkage.

\begin{figure}[htb!]
\centering
\includegraphics[width=2.0in]{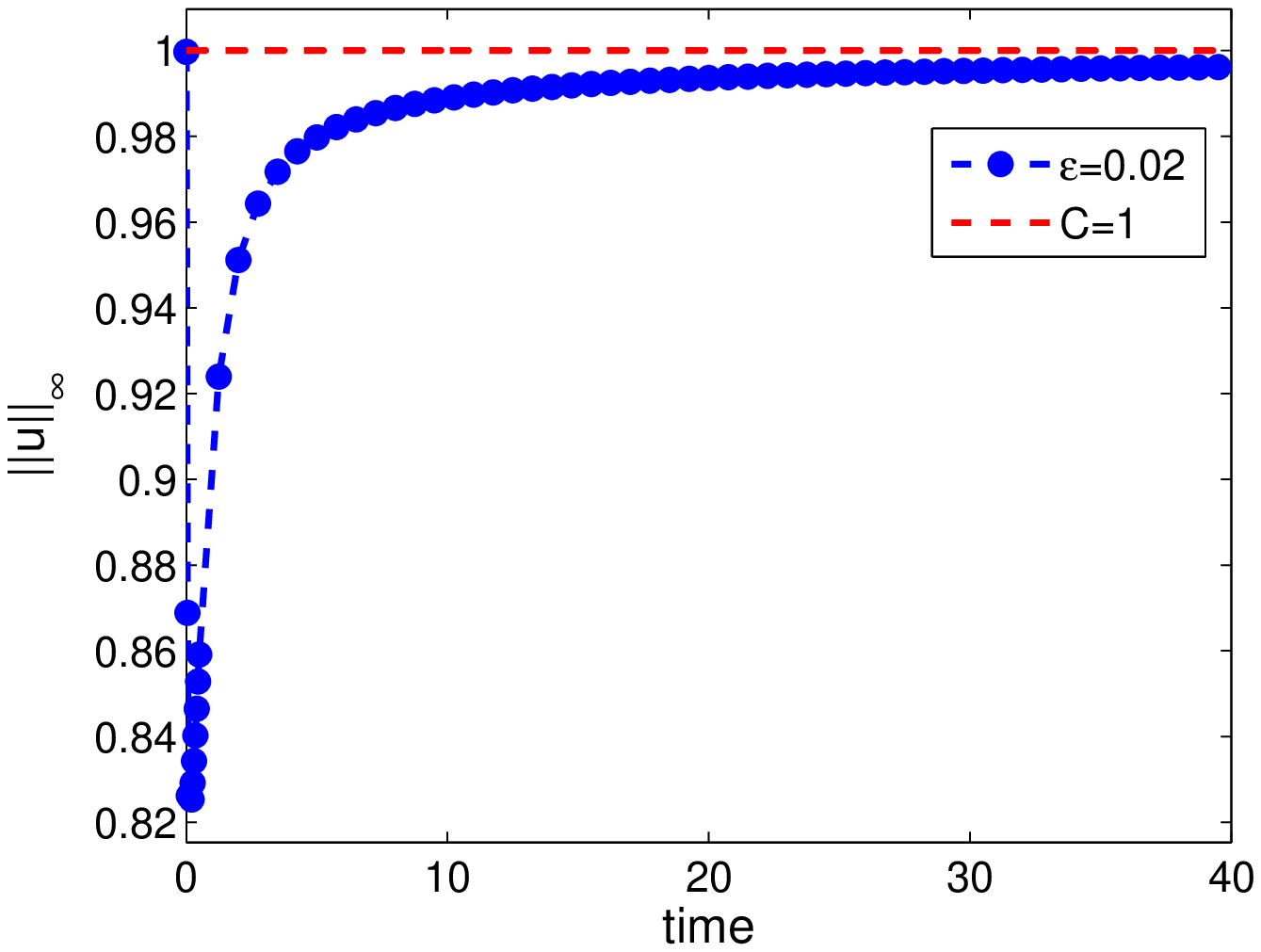}
\includegraphics[width=2.0in]{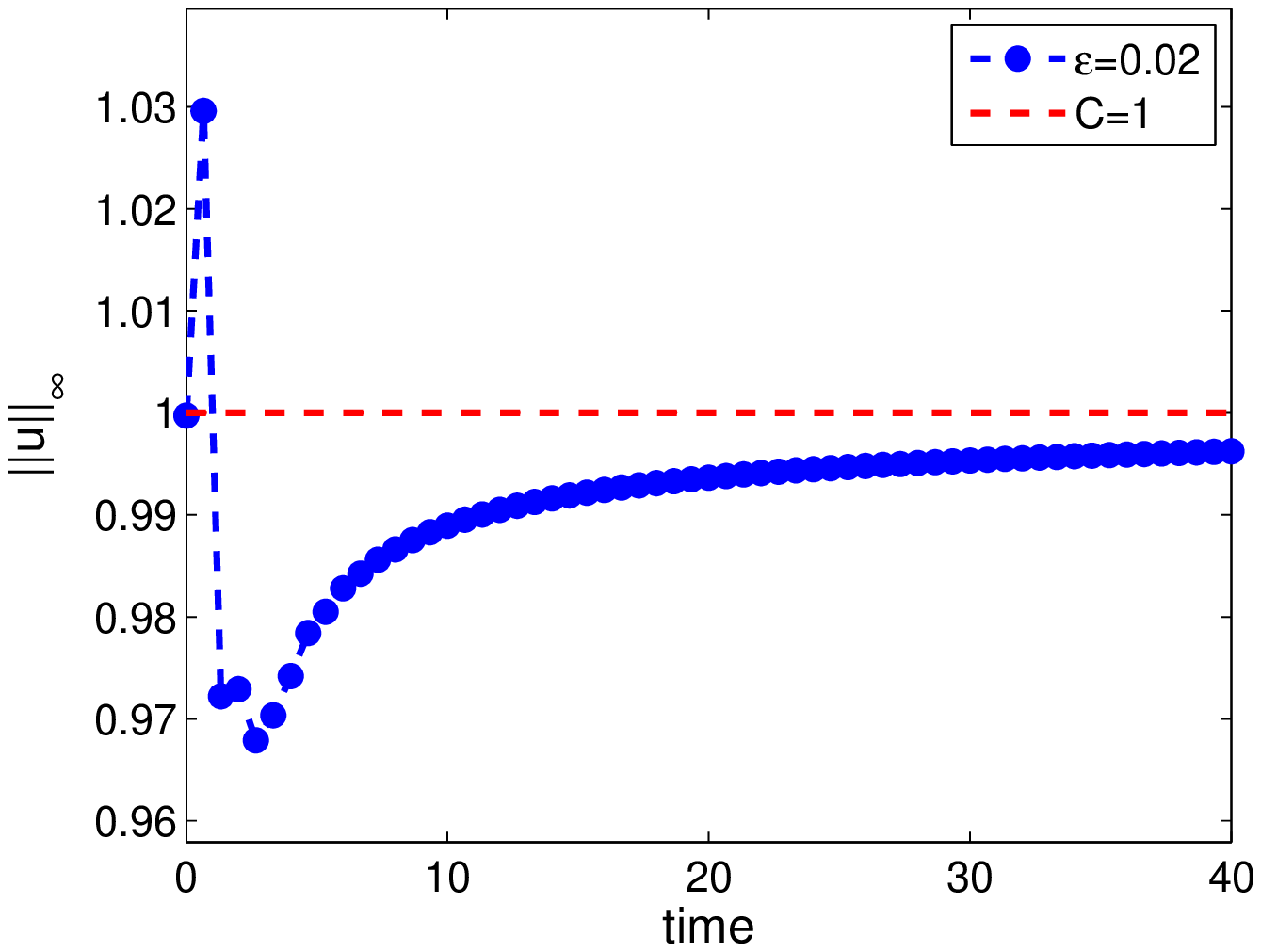}
\includegraphics[width=2.0in]{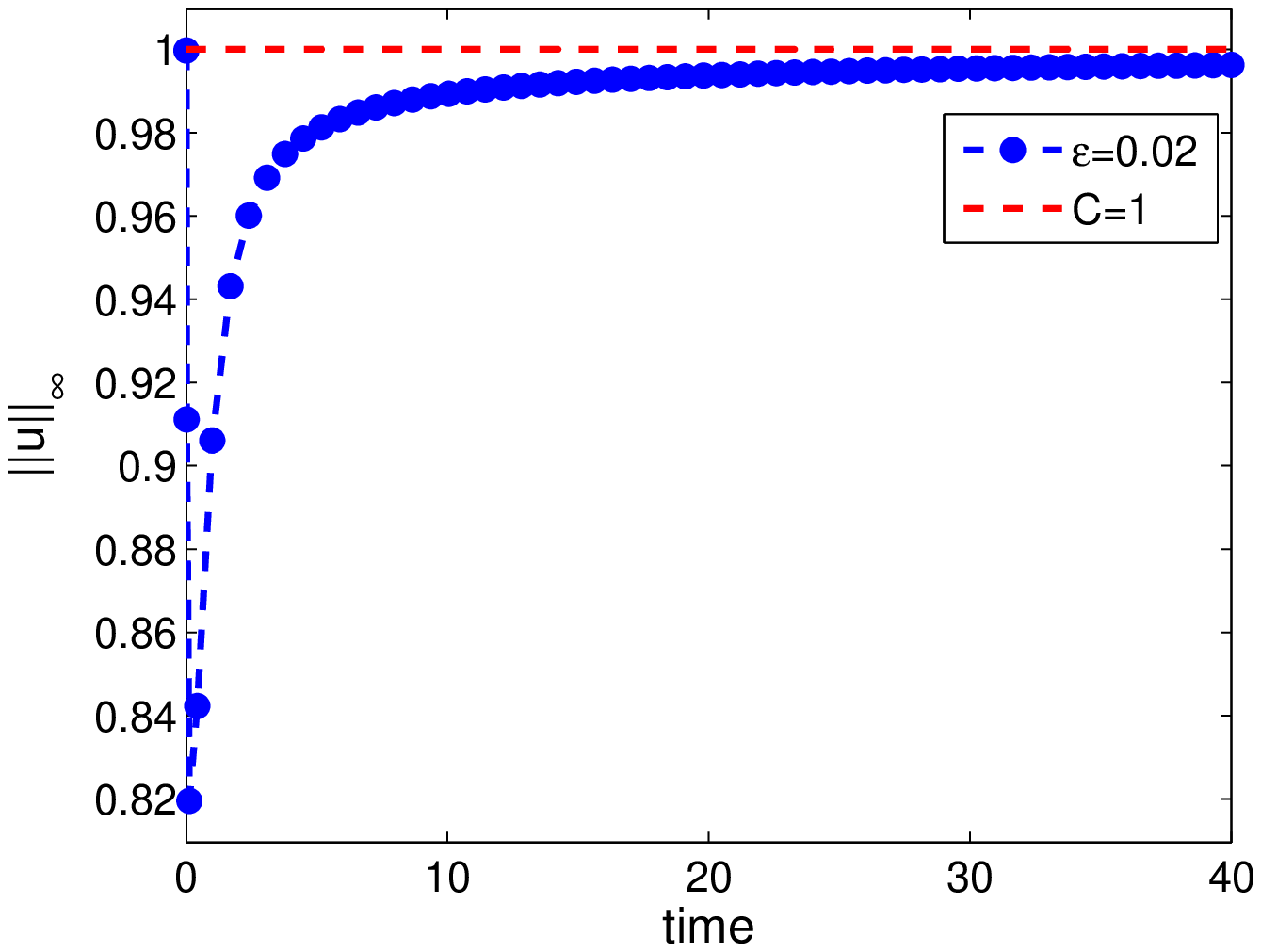}\\
\includegraphics[width=2.0in]{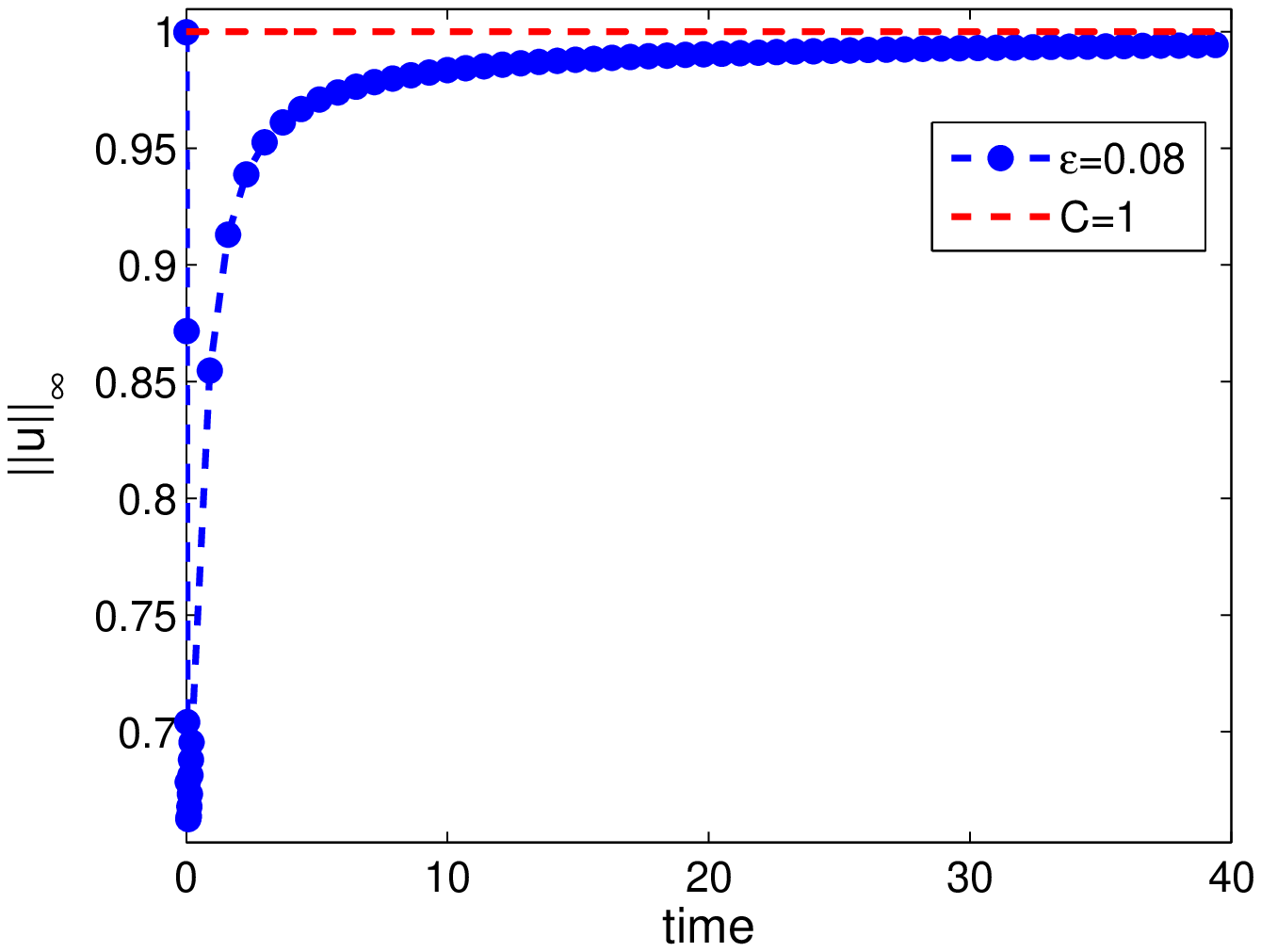}
\includegraphics[width=2.0in]{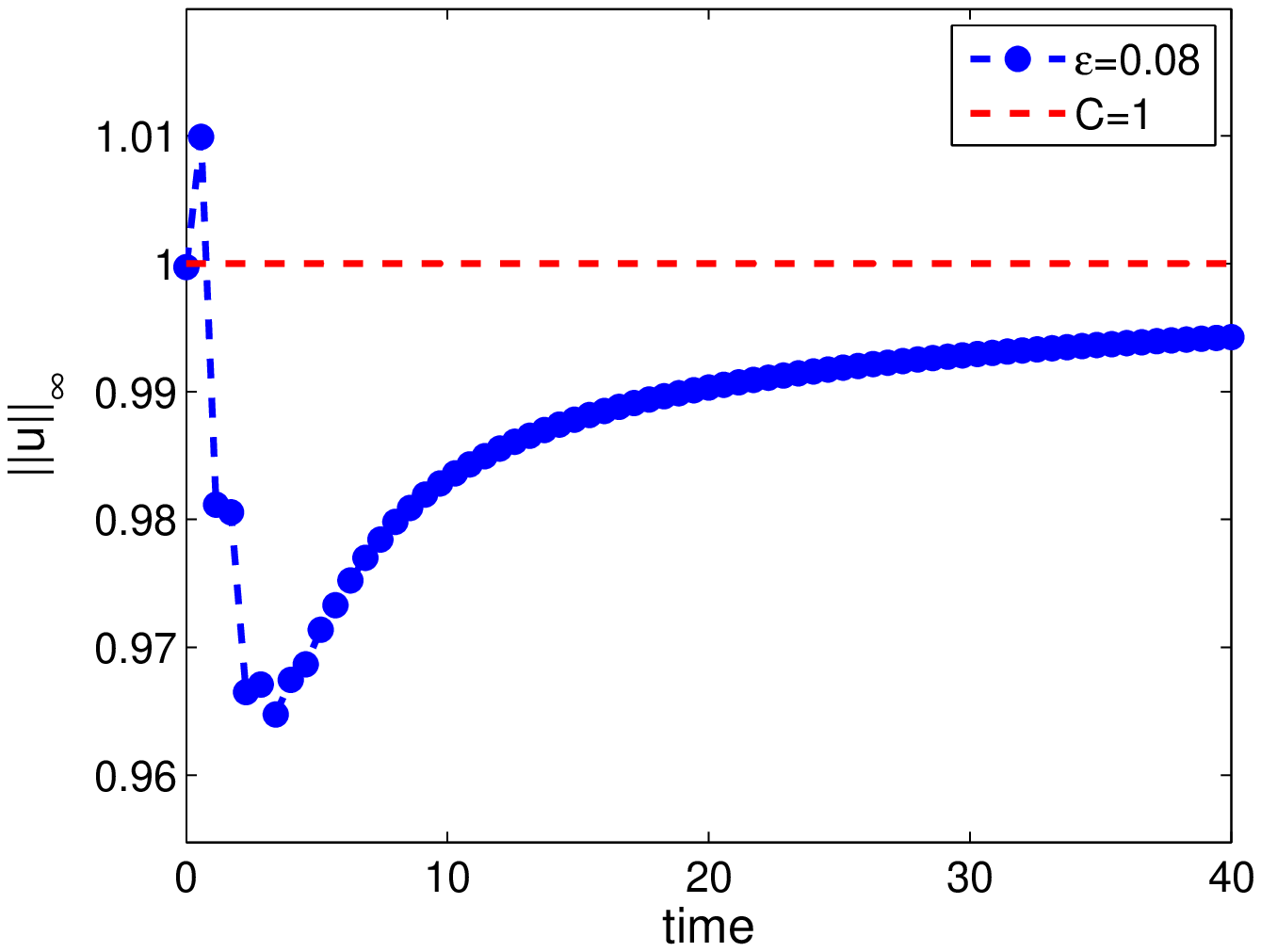}
\includegraphics[width=2.0in]{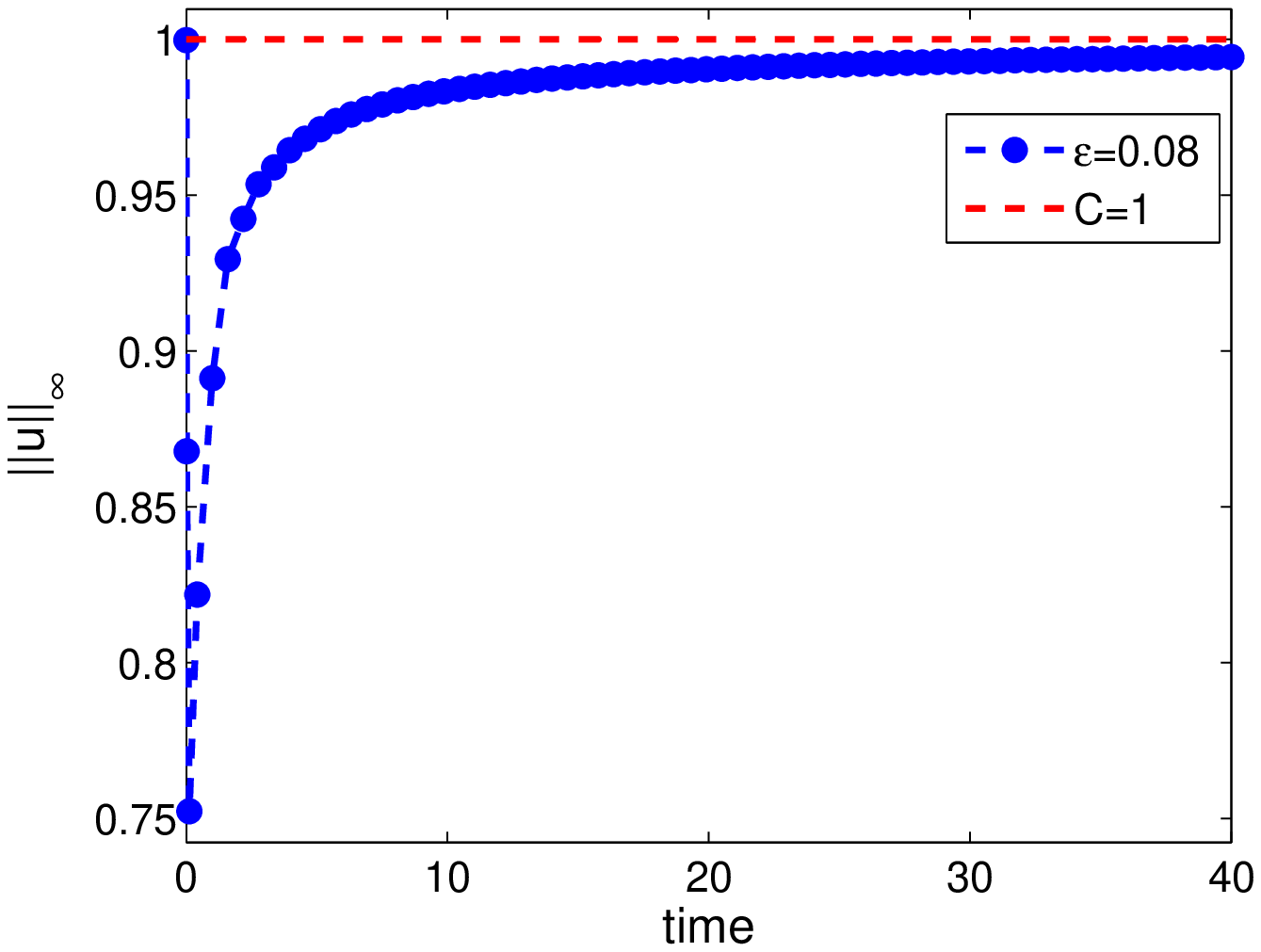}\\
\caption{The solution in maximum norm for $\varepsilon=0.02$ (top) and $0.08$ (bottom) on different time meshes.
Three meshes used in top figures: uniform mesh ($N=800$, $\tau=0.05$), uniform mesh ($N=60$, $\tau=0.67$) and graded mesh
($N=60$, $N_{0}=8$, $T_{0}=2$ and $\gamma=3$).
Three meshes adopted in bottom figures: uniform mesh ($N=2000$, $\tau=0.02$), uniform mesh ($N=70$, $\tau=0.57$) and graded mesh
($N=70$, $N_{0}=4$, $T_{0}=1$ and $\gamma=3$).}
\label{Maximum-Principle-Bound-Violate}
\end{figure}

\begin{example}\label{example:maximum principle}
We next consider $\partial_{t}^{\alpha}u=\varepsilon^{2}\Delta u-f(u)$ on $\Omega=(0,1)^{2}\times(0,40]$ with the fractional order $\alpha=0.7$.
The solution is computed with the spatial step $h=0.01$ using initial data
$u_{0}=0.95\times{rand}(\mathbf{x})+0.05$, where $rand(\cdot)$ generates a random number in $(0,1)$.
\end{example}
We use this example to examine the discrete maximum principle by two different diffusive coefficients $\varepsilon=0.02, 0.08$ and three different time-stepping approaches,
see Figure \ref{Maximum-Principle-Bound-Violate}. Notice that the graded meshes in the right figures
put $N_0$ grid points with a proper grading parameter $\gamma$ inside the starting cell $[0,T_0]$, cf. Example \ref{example:Accuracy-Test},
but use the uniform mesh with the time-step $\tau=(T-T_0)/(N-N_0)$ over the remainder interval $(T_0,T]$.
The try-and-error tests show that the maximum norm values are uniformly bounded by 1
provided the time-step size $\tau<0.67$ and $\tau<0.57$ for the two cases $\varepsilon=0.02$ and $\varepsilon=0.08$, respectively.
As seen, the time-step constraint \eqref{Time-Step-Constraint} is only sufficient to ensure the discrete maximum principle.

More interestingly, when the graded mesh is adopted near the initial time, the maximum norm values are still
bounded by 1 even for larger time-steps
($\tau=0.73$  for $\varepsilon=0.02$ in the top right figure in Figure \ref{Maximum-Principle-Bound-Violate})
in the remainder interval $(T_0,T]$. This shows that a good resolution of initial singularity
is also important to preserve the maximum principle.

\section{Conclusions}
We have proposed a second-order maximum principle preserving time-stepping scheme for the time-fractional Allen-Cahn equation under nonuniform time steps. Sharp maximum-norm error estimates the can reflect the temporal regularity are also presented. As our analysis is built on nonuniform time steps, we may resolve the intrinsic initial singularity by considering the graded meshes, and furthermore, we propose an adaptive time-stepping strategy for long-time simulations. Numerical experiments are presented to show the effectiveness of the proposed scheme.

We remark that the energy stability has not beed addressed in this work. Up to now we are unable to build up a discrete energy dissipation law for the second-order scheme {\eqref{Scheme-1}}-{\eqref{Scheme-2}}. As seen in \cite{Tang2018On}, the key issue is
to prove the positive semi-definite of the quadratic form $\sum_{k=1}^nw_k\sum_{j=1}^kA_{k-j}^{(k)}w_j$. In fact, we can show the energy stability under \textit{uniform} mesh using similar arguments as in \cite{Tang2018On}.
However, on a general \textit{nonuniform} mesh, it remains open to determine what kind of restrictions must be imposed on the discrete kernels $A_{n-k}^{(n)}$
so that $\sum_{k=1}^nw_k\sum_{j=1}^kA_{k-j}^{(k)}w_j$ is positive semi-definite. This will be part of our future studies.

\bibliographystyle{plain}
\bibliography{AllenChan_BibFile1,fractional,fractional_new}
\end{document}